\documentclass[12pt,a4paper,reqno]{amsart}
\usepackage{amsmath,amsthm,amscd}
\usepackage{amsfonts,amssymb}
\usepackage{tikz-cd}
\allowdisplaybreaks

\DeclareMathOperator{\mes}{mes}
\DeclareMathOperator{\esssup}{ess\,sup}

\newcommand{\be}{\begin{equation}}
\newcommand{\ef}{\end{equation}}
\chardef\bslash=`\\ 

\newtheorem*{thm*}{Theorem}

\theoremstyle{definition}

\newtheorem*{remark*}{Remarks}
\newtheorem*{defn*}{Definition}
\theoremstyle{remark}
\theoremstyle{remark*}

\newcommand{\wt}{\widetilde}
\newcommand{\wh}{\widehat}

 \renewcommand{\sectionmark}[1]{}
\renewcommand{\Im}{\operatorname{Im}}
\renewcommand{\Re}{\operatorname{Re}}
\newcommand{\iy}{\infty}

\newcommand{\qc} {quasiconformal}

\newcommand{\ve}{\varepsilon}

\newcommand{\fc}{\frac}
\newcommand{\Te} {Teichm\"{u}ller}
\newcommand{\const}{\operatorname{const}}

 \usepackage{amsfonts}
\newcommand{\field}[1]{\mathbb{#1}}

\newcommand{\dl}{\delta}
\newcommand{\D}{\field{D}}
\newcommand{\om}{\omega}
\newcommand{\z}{\zeta}
\newcommand{\ov}{\overline}
\newcommand{\vp}{\varphi}
\newcommand{\hC}{\widehat{\field{C}}}
\newcommand{\C}{\field{C}}
\newcommand{\R}{\field{R}}
\newcommand{\hR}{\widehat{\field{R}}}
\newcommand{\B}{\mathbf{B}}
\newcommand{\T}{\mathbf{T}}

\newcommand{\Hol}{\operatorname{Hol}}
\newcommand{\id}{\operatorname{id}}
\newcommand{\grad}{\operatorname{grad}}

\newcommand{\Belt} {\operatorname{Belt}}

\newcommand{\vk} {\varkappa}

\newcommand{\x} {\mathbf x}
\renewcommand{\a} {\alpha}
\newcommand{\ld}{\lambda}

\newcommand{\Teich}{\operatorname{Teich}}

\begin{document}

\title{The Grunsky operator and quasiconformality: old and new}
\author{Samuel L. Krushkal}

\begin{abstract} The Grunsky operator arises from univalence and plays a crucial role in geometric function theory. This operator also implies quasiconformal extendibility and has an intrinsic connection with Teichm\"{u}ller space theory and its interactions with complex analysis and pluripotential theory.

This paper surveys recent results in this field and simultaneously provides solutions of two old problems.

\end{abstract}

\date{\today\hskip4mm({GrOpQc(1).tex})}

\maketitle

\bigskip

{\small {\textbf {2020 Mathematics Subject Classification:} Primary: 30C55, 30C62, 30C75, 30F60, 32F45, 32G15; Secondary 30-08, 30F45, 31C10, 46G20}

\medskip

}

\bigskip

\markboth{Samuel L. Krushkal}{Grunsky operator and quasiconformality: old and new} \pagestyle{headings}

\bigskip\bigskip
\centerline{\bf 1. BACKGROUND}

\bigskip\noindent
{\bf 1.1.  The Grunsky operator on univalent functions and its generalizations}.
The classical Grunsky theorem \cite{Gr} states that a holomorphic function
$F(z) = z + \const + O(z^{-1})$ in a neighborhood of $z = \iy$ is continued to a univalent function
on the disk
$$
\D^* = \{z \in \hC = \C \cup \{\iy\}: \ |z| > 1\}
$$
if and only if the coefficients $\a_{m n}(f)$ of the expansion
 \be\label{1}
\log \fc{f(z) - f(\z)}{z - \z} =  \sum\limits_{m, n = 1}^\iy \a_{m
n} z^m \z^n \quad (z, \z \in \D^2)
\end{equation}
(where the principal branch of the logarithmic function is chosen) satisfy the inequality
 \be\label{2}
\Big\vert \sum\limits_{m,n=1}^\iy \ \sqrt{m n} \ \a_{m n}(f) x_m x_n \Big\vert \le 1
\end{equation}
for any sequence $\mathbf x = (x_n)$ from  the unit sphere $S(l^2)$ of the Hilbert space $l^2$ with
norm $\|\x\| = \bigl( \sum\limits_1^\iy |x_n|^2\bigr)^{1/2}$.

The class of all univalent $\hC$-holomorphic functions $f(z) = z + b_0 + b_1 z^{-1} + \dots$ on the
disk $\D^*$ is denoted by $\Sigma$.
The inversions $F_f(z) = 1/f(1/z)$ of zero free functions $f$ from $\Sigma$ form the canonical class
$S$ of univalent functions on the unit disk having the same Grunsky coefficients as $f$ (and accordingly the subclass $S_Q(\iy)$ with $f(\iy) = \iy$) and the Schwarzians $S_{F_f}$ run over a bounded
domain in the corresponding space $B(\D)$.
It is technically more convenient to deal with univalent functions on $\D^*$, but all results presented
below can be reformulated in terms of functions from $S$.

For technical convenience, all results  will be formulated for functions $f \in \Sigma_Q(0)$ and similar classes in arbitrary unbounded quasidisks; the corresponding functions $F_f$ on bounded quasidisks have similar features.

The technique of the Grunsky inequalities provides one of the most powerful methods in the classical geometric function theory, see, e.g., \cite{L}, \cite{Mi}, \cite{Po1}.

\bigskip
The aim of this paper is to investigate the interaction between the Grunsky operator and quasiconformality in connection with the crucial role of univalent functions with quasiconformal extension in complex analysis and Teichm\"{u}ller space theory. It also surveys the most topics in
this field.

The start point is that in the case of functions $f \in \Sigma$ admitting $k$-quasiconformal
extension to the whole sphere $\hC$, the inequality (2) is strengthened by
$$
\Big\vert \sum\limits_{m,n=1}^\iy \ \sqrt{m n} \ \a_{m n}(f) x_m x_n \Big\vert \le k;
$$
this was first established by K\"{u}hnau in \cite{Ku1}.

We shall consider the subclass $\Sigma_Q(0)$ from $\Sigma$ consisting of zero-free functions on
$\D^*$ which have quasiconformal extensions across the unit sphere $\mathbb S^1 = \{|z| = 1\}$
(hence to the whole Riemann sphere $\hC$) and complete this class in the topology of locally uniform convergence in $\D^*$ in the spherical metric on $\hC$.

The Beltrami coefficients of quasiconformal extensions are supported in the unit disk $\D =
\{|z| < 1\}$ and run over the unit ball
$$
\Belt(\D)_1 = \{\mu \in L_\iy(\C): \ \mu(z)|\D^* = 0, \ \ \|\mu\|_\iy  < 1\}.
$$
Each $\mu \in \Belt(\D)_1$ determines a unique homeomorphic solution to the Beltrami equation
$\ov \partial w = \mu \partial w$ on $\C$ (quasiconformal automorphism of $\hC$)  normalized by
$w^\mu(\iy) = \iy, \ (w^\mu)^\prime(\iy) = 1, \ w^\mu(0) = 0$, whose restriction to $\D^*$
belongs to $\Sigma_Q(0)$.

\bigskip\noindent
{\bf 1.2. Schwarzian derivatives and universal Teichm\"{u}ller space}.
Another quantity naturally associated with the univalent functions is the {\bf Schwarzian derivative}
$$
S_f(z) = \Bigl(\fc{f^{\prime\prime}(z)}{f^\prime(z)}\Bigr)^\prime
- \fc{1}{2} \Bigl(\fc{f^{\prime\prime}(z)}{f^\prime(z)}\Bigr)^2 \quad z \in \D^*.
$$
The Schwarzians $S_f$ of $f \in \Sigma$ belong to the complex Banach space $\B = \B(\D^*)$ of hyperbolically bounded holomorphic functions in the disk $\D^*$ with norm
$$
\|\vp\|_\B = \sup_{\D^*} \ (|z|^2 - 1)^2 |\vp(z)|
$$
and run over a  bounded domain in $\B$ modeling the {\bf universal Teichm\"{u}ller space} $\T$
(note that $\vp(z) = O(|z|^{-4})$ as $z \to \iy$).

The origin (the base point) $\vp = \mathbf 0$ of the space $\B$ corresponds to the identity map $f(z) \equiv z$.
This space is dual to the Bergman space $A_1(\D^*)$, a subspace of $L_1(\D^*)$ formed by integrable holomorphic functions (quadratic differentials $\vp(z) dz^2$ on $\D^*$.

The universal Teichm\"{u}ller space $\T = \Teich (\D)$ is the space of
quasisymmetric homeomorphisms of the unit circle $\mathbb S^1$ factorized by M\"{o}bius maps. Its topology and real geometry is determined by  Teichm\"{u}ller metric which naturally arises from extensions of those $h$ to the unit disk.

This space also admits the complex structure of a complex Banach manifold.
This structure is defined on $\T$ by factorization of the ball $\Belt(\D)_1$,
letting $\mu_1, \mu_2 \in \Belt(\D)_1$ be equivalent if the corresponding quasiconformal maps
$w^{\mu_1}, w^{\mu_2}$ coincide on the unit circle $\mathbb S^1 = \partial \D$
(hence, on $\ov{\D}$). Such $\mu$ and the corresponding maps $w^\mu$ are called $\T$-{\bf equivalent}.
The equivalence classes $[w^\mu]_\T$ are in one-to-one correspondence with the Schwarzian derivatives $S_{w^\mu}(z), \ \ z \in \D$.
The factorizing projection
$$
\phi_\T(\mu) = S_{w^\mu}: \ \Belt(\D^*)_1 \to \T
$$
is a holomorphic map from $L_\iy(\D^*)$ to $\B$. This map is a split submersion, which means that $\phi_\T$ has local holomorphic sections (see, e.g., \cite{EKK}, \cite{GL}, \cite{Le}).

The basic intrinsic metric on the space $\T$ is its {\bf Teichm\"{u}ller} metric
$$
\tau_\T (\phi_\T (\mu), \phi_\T (\nu)) = \frac{1}{2} \inf
\bigl\{ \log K \bigl( w^{\mu_*} \circ \bigl(w^{\nu_*} \bigr)^{-1} \bigr) : \
\mu_* \in \phi_\T(\mu), \nu_* \in \phi_\T(\nu) \bigr\};
$$
it is generated by the canonical Finsler structure on $\wt \T$ (in fact on the tangent bundle
$\mathcal T(\T) = \T \times \B$ of $\T$).

The {\bf Carath\'{e}odory} and {\bf Kobayashi} metrics on $\T$ are, as usually, the smallest
and the largest semi-metrics $d$ on $\T$, which are contracted by holomorphic
maps $h: \ \D \to \T$.
Denote these metrics by $c_\T$ and $d_\T$, respectively; then
$$
c_\T(\psi_1, \psi_2) = \sup \{d_\D(h(\psi_1), h(\psi_2)) : \
h \in \Hol(\T, \D)\},
$$
while $d_\T(\psi_1, \psi_2)$ is the largest pseudometric $d$ on $\T$ satisfying
$$
d(\psi_1, \psi_2) \le \inf \{d_\D(0, t) : \ h(0) = \psi_1, \ \text{and} \
h(t) = \psi_2 \ \ h \in \Hol(\D, \T)\},
$$
where $d_\D$ is the hyperbolic Poincar\'{e} metric on $\D$ of Gaussian curvature $- 4$.

By the Royden-Gardiner theorem,  the Kobayashi and Teichm\"{u}ller metrics are equal on all
Teichm\"{u}ller spaces (see, e.g., \cite{EKK}, \cite{GL}, \cite{Su}); so,
$$
k(f) = \tanh d_\T(\mathbf 0, S_f).
$$

\noindent
{\bf 1.3. Generalizations}.
The method of Grunsky inequalities was generalized in several directions, even to bordered Riemann surfaces $X$ with a finite number of boundary components, see \cite{L}, \cite{Mi}, \cite{SS}.
The case of arbitrary quasidisks from quasiconformal point of view is described in \cite{Kr9}.

Let $L \subset \C$ be a bounded (oriented) quasicircle separating the points $0$ and $\iy$ with the
interior and exterior domains $D$ and $D^*$ so that $0 \in D$ and $\iy \in D^*$.
Then the hyperbolic metric $\ld_D(z) |dz|$ of Gaussian curvature $-4$ on $D$ is estimated by
the Euclidean distance $e(z) = \inf_{\z \in L} |z - \z|$ by
$$
\fc{1}{4} \le \ld_D(z) e(z) \le 1.
$$
One naturally connects with these domains the corresponding collections of univalent functions: the class
$\Sigma(D^*)$ the collection of univalent functions on $D^*$ with hydrodynamical normalization
at $z = \iy$, i.e.,
$$
f(z) = b_0 + b_1 z^{-1} + \dots \quad \text{near} \ \ z = \iy,
$$
and the class univalent functions on $D$ normalized by $F(z) = z + a_2 z^2 + \dots$ near $z = 0$.
Let $\Sigma_k(D^*)$ and $S_k(D)$ be their subclasses formed by functions admitting $k$-quasiconformal extensions across the common boundary $L$ to the whole plane $\hC = \C \cup \{\iy\}$. The unions
$$
\Sigma_Q(D^*) = \bigcup_{k < 1} \Sigma_k(D^*), \quad S_Q(D) =  \bigcup_{k < 1} S_k(D)
$$
are dense in $\Sigma(D^*)$ and $S(D)$ in the topology of locally inform convergence on $D^*$ and $D$,
respectively. Without loss of generality, we can additionally normalize these extensions by $f^\mu(0) = 0$ and $F^\mu(\iy) = \iy$.

The Beltrami coefficients of these extensions run over the balls
$$
\Belt(D)_1 = \{\mu \in L_\iy(\C): \ \mu(z)|D^* = 0, \ \ \|\mu\|_\iy  < 1\}
$$
and $\Belt(D^*)_1 = \{\mu \in L_\iy(\C): \ \mu(z)|D = 0, \ \ \|\mu\|_\iy  < 1\}$.

The {\bf Grunsky-Milin} coefficients of functions $f \in \Sigma(D^*)$ are determined from expansion
 \be\label{3}
\log \fc{f(z) - f(\z)}{z - \z} = - \sum\limits_{m, n = 1}^\iy
\fc{\beta_{m n}}{\sqrt{m n} \ \chi(z)^m \ \chi(\z)^n},
\end{equation}
where $\chi$ denotes a conformal map of $D^*$ onto the disk $\D^*$ so that $\chi(\iy) = \iy, \ \chi^\prime(\iy) > 0$ , and the necessary and sufficient condition of univalence of a function
$f(z) = z + \const + O(1/z)$ near $z = \iy$ to continue to a univalent function on $D^*$ is
$$
\Big\vert \sum\limits_{m,n = 1}^{\iy} \ \beta_{mn} \ x_m x_n \Big\vert \le 1
$$
for any sequence ${\mathbf x} = (x_n)$ from the unit sphere $S(l^2) = \{\x \in l^2: \ \|\x\|= 1\}$
(cf. \cite{Gr}, \cite{Kr2}, \cite{Mi}, \cite{Po1}). The quantity
$$
\vk_{D^*}(f) = \sup \Big\{ \Big\vert \sum\limits_{m,n = 1}^{\iy} \ \beta_{mn} \ x_m x_n \Big\vert : \
{\mathbf x} = (x_n) \in S(l^2)\Big\}
$$
is called the {Grunsky norm} of $f$  and is majorated by the {\bf Teichm\"{u}ller norm} $k(f)$
equal to the smallest dilatation $k(f^{\mu_0})$ among quasiconformal extensions of function $f$ onto $D$;
moreover, on the open dense subset of $\Sigma_Q(D^*)$ the strict inequality $\vk_{D^*}(f) < k(f)$ is valid.

For any quasidisk $D^*$, the set of $f$ with $\vk_{D^*}(f) = k(f)$ is sparse, but such functions play a crucial role in many applications.

In the case $D^* = \D^*$, we have the canonical Grunsky norm, which will be denoted by $\vk(f)$.

\bigskip
The generalized Grunsky coefficients $\beta_{m n} (f^\mu)$ of $f^\mu \in \Sigma(D^*)$ determine for every fixed $\x = (x_n) \in l^2$ with $\|\x\| = 1$ the holomorphic map
 \be\label{4}
h_{\x}(\mu) =
\sum\limits_{m,n=1}^\iy \ \beta_{m n} (f^\mu) x_m x_n : \ \Belt(D)_1 \to \D,
\end{equation}
and $\sup_{\x} |h_{\x}(f^\mu)| = \vk_{D^*}(f^\mu)$.
The holomorphy of functions (4) follows from the holomorphy of coefficients $\beta_{m n}$ with
respect to Beltrami coefficients $\mu \in \Belt(D)_1$ using the estimate
 \be\label{5}
\Big\vert \sum\limits_{m=j}^M \sum\limits_{n=l}^N \ \beta_{mn} x_m x_n \Big\vert^2
\le \sum\limits_{m=j}^M |x_m|^2 \sum\limits_{n=l}^N |x_n|^2
\end{equation}
which holds for any finite $M, N$ and $1 \le j \le M, \ 1 \le l \le N$.
This estimate is a simple corollary of Milin's univalence theorem (cf. \cite{Mi}, p. 193;
\cite{Po1}, p. 61).

Both norms $\vk(f)$ and $k(f)$ are continuous plurisubharmonic functions of $S_f$ on the space $\T$
(in $\B$ norm); see, e.g. \cite{Kr9}.

In particular, for univalent functions $f(z) = z + a_2 z^2 + \dots$ in an arbitrary disk
$\D_r^* = \{z \in \hC: \ |z| > r\}, \ 0 < r < \iy$, the corresponding function (1) on $(z, \z) \in \D_r^2$
provides the Grunsky norm
$$
\vk_r(f) = \sup \Big\{\Big\vert \sum\limits_{m,n=1}^\iy \ \sqrt{mn} \ \a_{m n} r^{m+n} x_m x_n \Big\vert: \ \mathbf x = (x_n) \in S(l^2) \Big\},
$$
and accordingly, the holomorphic maps
  \be\label{6}
 h_{\x,r}(S_f) = \sum\limits_{m,n=1}^\iy \ \sqrt{m n} \ \a_{m n} (S_f) r^{m+n} \ x_m x_n : \ \T \to \D
 \end{equation}
with
 \be\label{7}
\sup_{\x\in S(l^2)} h_{\x,r}(S_f) = \vk_r(f).
\end{equation}

Note also that the Grunsky (matrix) operator
$$
\mathcal G(f) = (\sqrt{m n} \ \a_{mn}(f))_{m,n=1}^\iy
$$
acts as a linear operator $l^2 \to l^2$ contracting the norms of elements $\x \in l^2$;
the norm of this operator equals $\vk(f)$ (cf. \cite{G1}).

\bigskip\noindent
{\bf 1.4. The root transform}.
One can apply to $f \in \Sigma_Q(\iy)$ (and to its inversion $F_f$) the rotational conjugation
$$
\mathcal R_p: \ f(z) \mapsto f_p(z) := f(z^p)^{1/p} = z + \fc{a_2}{p} z^{p-1} + \dots
$$
with integer $p \ge 2$, which transforms $f \in S(\iy)$ into $p$-symmetric univalent functions accordingly
to the commutative diagram
$$
\begin{tikzcd}
\wt \C_p \arrow{d}{\pi_p}  \arrow{r}{\mathcal R_p f}
         &\wt \C_p \arrow{d}{\pi_p}      \\
\hC \arrow{r}{f}
         &\hC
\end{tikzcd}
$$
where $\wt \C_p$ denotes the $p$-sheeted sphere $\hC$ branched over $0$ and $\iy$, and the projection $\pi_p(z) = z^p$.

This transform acts on $\mu \in \Belt(\D)_1$ and $\psi \in L_1(\D)$ by
 \be\label{8}
\mathcal R_p^* \mu = \mu(z^p) \ov z^{p-1}/z^{p-1}, \quad \mathcal
R_p^* \psi = \psi(z^p) p^2 z^{2p-2};
\end{equation}
then
$$
k(\mathcal R_p f) = k(f), \quad  \vk(\mathcal R_p f) \ge \vk(f).
$$
It was established by Grinshpan in \cite{G1} that the root transform does not decrease the
Grunsky norm, i.e.
$$
\vk_p(f) := \vk(f_p) \ge \vk(f)
$$
(this also follows from the K\"{u}hnau-Schiffer theorem on reciprocity of the Grunsky norm to
the least positive Fredholm eigenvalue of the curve $L = f(|z| = 1)$; see \cite{Ku2}, \cite{Sc}).

Note that the sequence $\vk_p(f), \ p = 2, 3, \dots$, is not necessarily nondecreasing. For example, for the function $F(z) = z/(1 + t z)^2$ with $|t| \le 1$, we have
$$
F_p(t) = z/(1 + t z^p)^{2/p},
$$
and $\vk(F_p) = |t| = k(f_p)$ for even $p$, while $\vk(F_p)< |t| = k(F_p)$ for any odd $p \ge 3$ (see, e.g. \cite{Kr4}).

\bigskip
The Grunsky coefficients of every function $\mathcal R_p f$ also are polynomials of
$b_1, \dots \ , b_l$, which implies, similar to (4), holomorphy of maps
 \be\label{9}
h_{\x,p}(\mu) = \sum\limits_{m,n=1}^\iy \ \sqrt{m n} \ \a_{m n} (\mathcal R_p f^\mu) \  x_m x_n : \
\Belt(D)_1 \to \D
\end{equation}
for any fixed $p$ and any $\x = (x_n) \in S(l^2)$, and
$\sup_{\x\in S(l^2)} h_{\x,p}(\mu) = \vk(\mathcal R_p f^\mu)$.

Every function $h_{\x,p}(\mu)$ descends to a holomorphic functions on the space $\T$, which implies that all Grunsky norms $\vk(\mathcal R_p f^\mu)$ are continuous and plurisubharmonic on $\T$ .

\bigskip\bigskip
\centerline{\bf 2. BASIC FEATURES OF GRUNSKY NORM}

\bigskip\noindent
{\bf 2.1}.
The following important general theorem on restoration of Grunsky norm by abelian holomorphic differentials was established recently in \cite{Kr18}. This theorem is underlying for many results
established below.

Let $D$ and $D^*$ be the interior and exterior domains of a bounded quasicircle $L$.

\bigskip\noindent
{\bf Theorem 1}. {\it For any univalent function $f = f^{t \mu_0}$ on $D^*$, we have the equality
 \be\label{10}
\vk(f^{t \mu_0}) = |t| \fc{|t| + \a(f)}{1 + \a(f) |t|} = \a(f) |t| + (1 - \a(f)^2) |t|^2 + \dots ,
\end{equation}
where $\mu_0$ is an extremal Beltrami coefficient among quasiconformal extensions $f^\mu$ of $f$ onto $D$, and
 \be\label{11}
\a(f) = \sup \ \Big\{ \Big\vert \iint_D \mu_0(z) \psi(z) dx dy \Big\vert: \ \psi \in A_1^2(D), \ \|\psi\|_{A_1} = 1 \Big\} \quad (z = x + iy),
\end{equation}
$A_1(D)$ denotes the subspace in $L_1(D)$ formed by integrable holomorphic functions
(quadratic differentials) on $\D$, and $A_1^2(\D)$ is its subset consisting
of $\psi$ with zeros even order in $D$, i.e., of the squares of holomorphic functions. }

\bigskip
The proof of this theorem is very complicated. It involves Golusin's improvement of Schwarz's lemma,
some deep Milin's estimates for kernels and orthogonal systems in multiply connected domains and estimating the metric $\ld_\vk$ of generalized Gaussian curvature $\kappa_{\ld_\vk} \le - 4$ by
Minda's maximum principle (the details see in \cite{Kr17}).

If $D^*$ is the canonical disk $D^* = \{|z| > 1\}$, then due to \cite{Kr4}, every $\psi \in A_1^2(D)$ has the form
$$
\psi(z) = \fc{1}{\pi} \sum\limits_{m+n = 2}^\iy \sqrt{m n} \ x_m x_n z^{(m+n)}
$$
and $\|\psi\|_{A_1(\D)} = \|\x\|_{l^2} = 1, \ \x = (x_n)$.

In this case the quantity $\a(f)$ for every $f \in \Sigma_Q$ is represented in the form
$$
\a(f) = \sup_{\x=(x_n)\in S(l^2)} \ \fc{1}{\pi \|\mu_0\|_\iy} \Big\vert
\iint\limits_{|z|<1} \mu_0(z) \ \sum\limits_{m+n\ge 2}^\iy \sqrt{m n} \ x_m x_n z^{m+n-2} dx dy \Big\vert,
$$
where $\mu_0$ is any extremal Beltrami coefficient for $f$ (cf. \cite{Kr4}).

\bigskip
In contrast to (10), a Beltrami coefficients $\mu_0 \in \Belt(D^*)_1$ is {\bf extremal} in its equivalence class (has the minimal $L_\iy$-norm) if and only if
$$
\|\mu_0\|_\iy = \sup \ \Big\{ \Big\vert \iint\limits_{D^*} \mu(z) \psi(z) dx dy \Big\vert: \ \psi \in A_1(D), \ \|\psi\|_{A_1} = 1\Big\}.
$$
This result is given by the Hamilton-Krushkal-Reich-Strebel theorem and is underlying for many
results in quasiconformal anlysis and Teichm\"{u}ller space theory.

\bigskip
Together with (10), this theorem yields the following important consequence (cf. \cite{Kr8}, \cite{Kr11}, \cite{Ku3}, \cite{St1}), which also is crucial in many applications.

\bigskip\noindent
{\bf Corollary 1}. {\it The equality $\vk(f) = k(f)$ is valid only when the both quantities
$k(f)$ and $\vk_{D^*}(f)$ coincide with $\a(f)$ given by (11).
Moreover, if $\vk_{D^*}(f) = k(f)$ and the equivalence class of $f$ is a {\bf Strebel point} of $\T$,
which means that this class contains the Teichm\"{u}ller extremal extension  $f^{k|\psi_0|/\psi_0}$
with $\psi_0 \in A_1(D)$, then necessarily $\psi_0 = \om ^2 \in A_1^2(D)$.
}

\bigskip
An important fact is that the Strebel points are dense in any Teichm\"{u}ller space, see \cite{GL}.

\bigskip\noindent
{\bf 2.2}.
We also mention the following deep theorem of Pommerenke and Zhuravlev given in \cite{Po1}; \cite{KK1}, pp. 82-84; \cite{Zh1}.

\bigskip\noindent
{\bf Theorem 2}. {\it Any univalent function $f \in S$ with $\vk(f) \le k < 1$ has $k_1$-quasiconformal extensions to $\hC$ with}  $k_1 = k_1(k) \ge k$.

\bigskip
Some explicit estimates for $k_1(k)$ are established, for example, in \cite{Kr8}, \cite{Ku4}.

Theorem 1 and one of results established below provide explicitly the minimal value of the
possible bounds $k_1(k)$ (see Section {\bf 5.4}).

\bigskip\noindent
{\bf 2.3. Complex homotopy}. One can define for each $f \in \Sigma(0)$ the complex homotopy
A close operation is determined the holomorphic homotopy
$$
f_t(z) = t f(t^{-1} z): z + t b_0 + t^2 b_1 z^{-1} + \dots: \ \D^* \times \D \to \C
$$
with $|t| \le 1$, which determines for $|t| < 1$ a holomorphic map of the space $\T$ into itself by
$$
S_f(z) = S_{f_t}(z) = t^2 S_f(z/t).
$$
The latter is obtained by applying, for example, the following lemma of Earle \cite{Ea}.

\bigskip\noindent
{\bf Lemma 1}. {\it Let $E, T$ be open subsets of complex Banach spaces $X, Y$ and $B(E)$ be a Banach space of holomorphic functions on $E$ with sup norm. If $\phi(x, t)$ is a bounded map $E \times T \to B(E)$ such that $t \mapsto \phi(x, t)$ is holomorphic for each $x \in E$, then the map $\phi$ is holomorphic.}

Note that both quantities $k_{f_t}$ and $\vk(f_t)$ are circularly symmetric logarithmically subharmonic functions of
$t \in \D$; hence, these functions are continuous in $|t|$ on $[0, 1]$ and convex with respect to $\log |t|$, and
$$
\lim\limits_{|t| \to 1} k(f_t) = k(f), \quad \lim\limits_{|t| \to 1} \vk(f_t) = \vk(f).
$$

\bigskip\noindent
{\bf 2.4}.
Combining Theorem 1 with the K\"{u}hnau-Schiffer theorem mentioned above, one obtains explicitly the first Fredholm eigenvalue $\rho_L$ of every quasiconformal curve $L \subset \hC$. Namely, we have

\bigskip\noindent
{\bf Theorem 3}. {\it  For any quasicircle $L \subset \hC$,
 \be\label{12}
\fc{1}{\rho_L} = \sup_{\psi \in A_1^2(\D),\|\psi\|_{A_1} =1} |t|
\fc{|t| + |\langle \mu_0, \psi\rangle_\D|}{1 + |\langle \mu_0, \psi\rangle_\D| |t|},
\end{equation}
where $\mu_0$ is an extremal Beltrami coefficients of the appropriately normalized exterior conformal mapping function $f$ of $\D^*$ onto the exterior domain of the curve $L$ on which the Teichm\"{u}ller norm of $f$ is attained,
and $\langle \mu_0, \psi\rangle_\D = \iint_\D \mu_0 \psi dx dy$.

Moreover, letting in (10) $\mu_0 = t |\psi|/\psi$ with $|t| < 1$ and $\psi \in A_1(\D)$, one obtains a dense subset of possible eigenvalues $\rho_L$. }

\bigskip
If the equivalence class of $f$ is a Strebel point and $\vk(f) = k(f)$, then both sides of (10)
are equal to the reflection coefficient of the curve $L$.

A significant feature of Theorem 3 is that it is valid for all quasicircles, whereas all previous
results on Fredholm eigenvalues have been obtained only for much narrower collection of curves associated
with functions having equal Grunsky and Teichm\"{u}ller norms.

\bigskip\noindent
{\bf 2.5. Example 1}. To illustrate Theorem 1, consider the ellipse $\mathcal E$ with the foci at $-1, 1$ and semiaxes $a, b \ (a > b)$, and denote its exterior by$\mathcal E^*$.
An orthonormal basis in the Hilbert space $A_2(\mathcal E)$ of the square integrable holomorphic functions on $D$, is formed by the polynomials
$$
P_n(z) = 2 \sqrt{\fc{n + 1}{\pi}} \ (r^{n+1} - r^{-n-1}) \ U_n(z),
$$
where $r = (a + b)^2$ and $U_n(z)$ are the Chebyshev polynomials of the second kind,
$$
U_n(z) = \fc{1}{\sqrt{1 - z^2}} \ \sin [(n + 1) \arccos z], \quad n = 0, 1, \dots,
$$
(see \cite{Ne}), and one obtains by the Riesz-Fisher theorem that every function
$\om \in A_2(\mathcal E)$ is of the form
$$
\om(z) = \sum\limits_0^\iy x_n P_n(z), \quad \x = (x_n) \in l^2,
$$
with $\|\psi\|_{\mathcal E}^2 = \|\x\|_{l^2}$. Then Theorem 1 implies that for any $f \in \Sigma_{\mathcal E^*}$ its Grunsky norm
$$
\vk_{\mathcal E^*}(f) = \sup_{\|\om\|_{A^2(\mathcal E) = 1}} |t| \fc{|t| + \a(f)}{1 + \a(f) |t|}.
$$

\newpage
\centerline{\bf 3. RESTORATION OF QUASICONFORMAL DILATATION}
\centerline{\bf BY GRUNSKY NORM}

\bigskip\noindent
{\bf 3.1. Grinshpan's conjecture}. Theorem 1 shows that Grunsky norm of univalent functions
is completely determined by their extremal quasiconformal extensions.

It is naturally to consider the inverse question:

{\it In which extent the Grunsky norm delimits the dilatation (Teichm\"{u}ller norm) of univalent functions ?}

As was mentioned, generically $\vk_{D^*}(f) < k(f)$. The question is whether the Teichm\"{u}ller norm also is determined by original conformal map of $D^*$ (compare with Ahlfors' problem presented below
in Section {\bf 3.3}).

We provide the answer to this question in the case, when $D^*$ is the canonical disk $\D^*$.

The following deep conjecture was posed by A. Grinshpan in \cite{G2} (in an equivalent form):

\bigskip\noindent
{\it For every function $f \in \Sigma(0)$, we have the equality}
 \be\label{13}
\limsup\limits_{p\to \iy} \vk_p(f) = k(f).
\end{equation}
This conjecture arose in the theory of univalent functions with quasiconformal extensions and relates
to important problems of geometric complex analysis and of Teichm\"{u}ller space theory, because
geometrically the equality (13) means the equality of the Carath\'{e}odory and Kobayashi metrics on universal Teichm\"{u}ller space.

A slightly modified version of this conjecture was proven in \cite{Kr10}. We simplify here
this proof. Note that this version has the same geometric meaning.

\bigskip
Let us first mention that one cannot replace in the statement of Theorem 1 the assumption $f(z) \in \Sigma(0)$ by $f(z) \in \Sigma_Q$, i.e., omit the third normalization condition, because without this condition the root transform $\mathcal R_p$ can increase the Teichm\"{u}ller norm.

As an {\bf example}, consider the extremal map $g_r$ in Teichm\"{u}ller's
Verschiebungssatz with minimal dilatation among quasiconformal automorphisms of the unit disk, which are identical on the boundary circle and move the origin into a given point
$- r \in (-1, 0)$. Its Beltrami coefficient $\mu_0 = k |\psi_0|/\psi_0$ is defined by $\psi_0$, which is holomorphic and does not vanish on $\D \setminus \{0\}$ and has simple pole at $0$. This $\psi_0$ is orthogonal to all holomorphic quadratic differentials on $\D$ with respect to pairing
$$
\langle \vp, \psi\rangle = \iint_\D (1 - |z|^2)^2 \vp(z) \ov{\psi(z)} dx dy.
$$
For small $r$, the dilatation $k(g_r) = r/2 + O(r^2)$ (the corresponding formula in \cite{Te}, p. 343,
for extremal dilatation contains an error).

This map $g_r$ extends trivially to a quasiconformal map of $\hC$ by $g_r(z) = z$ for $|z| \ge 1$.
Consider the translated map $f_r(z) = g_r(z) + r$. Its restriction to $\D^*$ has dilatation
$k(f_r) = k(g_r) = 0$.

In contrast, the dilatation of the squared map $f_{r,2} := \mathcal R_2 f_r$ equals $r$, since
the differential $\mathcal R_2^* \psi_0$ is holomorphic on $\D$ and therefore the Beltrami  coefficient $\mathcal R_2^* \mu_0 = \mu_{f_{r,2}}$ is extremal for the boundary values $f_{r,2}|\mathbf S^1$ (note that $f_{r,2}(z) = z + r/(2 z) + \dots$ for $|z| > 1$), and
$$
\vk(f_{r,2}) = k(f_{r,2}) = r + O(r^2), \ r \to 0.
$$
Thus,
$$
\limsup\limits_{p\to \iy} \vk_p (f_r) = \vk_2 (f_r) > k(f_r).
$$
Another example was constructed by K\"{u}hnau.

\bigskip\noindent
{\bf 3.2. General restoration theorems}. We shall use the truncated Beltrami coefficients
$$
\mu_\rho(z) = \begin{cases}  \mu(\rho z),  \ \ &|z| > 1,   \\
                             0,     & |z| < 1,
\end{cases}
$$
with fixed $\rho, \ 0 < \rho < 1$;  the map $\mu \mapsto \mu_\rho$ generates a linear self-map of $\Belt(\D)_1$.

The answer to the problem stated above is given by the following theorem having many important consequences.

\bigskip\noindent
{\bf Theorem 4}. {\it Every univalent function $f(z) \in \Sigma(0)$ with Grunsky norm $\vk(f) < 1$ admits quasiconformal extension $f^\mu$ to the whole sphere $\hC$ with dilatation
 \be\label{14}
k \ge \wh \vk(f) := \sup_{\mu\in [\sigma_a \circ f]}
\ \limsup\limits_{r\to 1} \ \sup_p \sup_{\psi\in A_1^2(\D^*),
\|\psi\|_{A_1}=1} \
 \Big\vert \iint_{\D^*} \mathcal R_p^*\mu_r(z) \psi(z) dx dy \Big\vert
\end{equation}
This lower admissible bound is sharp and in the sense that it cannot be replaced by a smaller quantity for every $f$.
}

\bigskip
The quantity $ \wh \vk(f)$ can be regarded as the \textbf{outer limit Grunsky norm}
of $f$ on $\D$.

To clarify (14), note that if $\mu$ is of Teichm\"{u}ller form, $\mu = k |\psi|/\psi$, its root transform generates by (8) the Teichm\"{u}ller map $\mathcal R_p^* f= f^{k |\mathcal R_p^* \psi|/\mathcal R_p^* \psi}$ determined by quadratic differential
$\mathcal R_p^* \psi(z) = \psi(z^p) p^2 z^{2p-2}$.
So if $\psi$ has zero at a point $z_0 \in \D \setminus \{0\}$, then  $\mathcal R_p^* \psi(z) = 0$ at
the points $z = z_0^{1/p}$, and $|z_0^{1/p}| \to 1$ as $p \to \iy$.

Note also that the limit $\wh \vk(f)$ as a function of $S_f$ followed by its upper semicontinuous normalization is plurisubharmonic on the space $\T$.

\bigskip
Now we outline the proof of Theorem 4. First assume that $\mu_0$ is of Teichm\"{u}ller type.
Fix $\rho_j$ arbitrarily close to $1$ and pick appropriate large $p_j$ so that all zeros of odd order
of $\mathcal R_{p_j}^* \psi_0$ are placed in the annulus $\{1 < |z| < 1/\rho_j\}$.
Then, taking the truncated Beltrami coefficients $(\mathcal R_{p_j}^* \mu_0)_{\rho_j}$ for
$$
\mathcal R_{p_j}^* \mu_0 = k |\mathcal R_{p_j}^* \psi_0|/\mathcal R_{p_j}^* \psi_0,
$$
vanishing in the disks $\D_{1/\rho_j}$, one obtains by applying Corollary 1 the relation
$$
\vk(f^{(\mathcal R_{p_j}^* \mu_0)_{\rho_j}}) =
\sup_{(x_n) \in S(l^2)} \Big\vert \sum\limits_{m,n=1}^\iy \ \sqrt{m n} \
\a_{m n}(f^{(\mathcal R_{p_j}^* \mu_0)_{\rho_j}}) \  \rho_j^{m+n} x_m x_n \Big\vert.
$$

Using this equality, one can find the appropriate sequences $\{\rho_n\} \to 1, \ \{p_n\} \to \iy$
and $\{\psi_n\} \in A_1^2(\D^*)$ with $\|\psi_n\|_{A_1} = 1$ such that in the limit as
$n \to \iy$ the above relations result in the equality
$$
k(f^{\mu_0}) = \wt \vk(f^{\mu_0}) = \limsup\limits_{\rho\to 1} \ \sup_p \sup_{\psi\in A_1^2(\D^*), \|\psi\|_{A_1}=1} \
 \Big\vert \iint_{\D^*} \mathcal R_p^*\mu_{0 \rho}(z) \psi(z) dx dy \Big\vert,
$$
which proves Theorem 4 for functions with Teichm\"{u}ller extensions.

In the case of the generic $f(z) \in \Sigma(0)$, one can apply the above argument to their homotopy functions $f_r(z) = r f(z/r)$ with $r < 1$. The properties of norms $\vk_p(f_r)$ and $k(f_r)$ indicated in  Section {\bf 2.3} allows one to select their appropriate subsequences $\{\rho_{n_j}\}, \ \{p_{n_j}\}, \ \{\psi_{n_j}\}$ so that
 \be\label{15}
k(f) = \wt \vk(f) = \lim\limits_{\rho_{n_j}\to 1} \ \sup_{p_{n_j}}
\sup_{\psi_{n_j} \in A_1^2(\D^*), \|\psi_{n_j}\|_{A_1}=1} \
 \Big\vert \iint_{\D^*} \mathcal R_p^*\mu_{0 \rho_{n_j}}(z) \psi_{n_j}(z) dx dy \Big\vert,
\end{equation}
which implies the relations (14) for any $f \in \Sigma(0)$.

\bigskip\noindent
{\bf 3.3. Asymptotically conformal functions}.
Theorem 4 provides that the dilatation $k(f)$ is attained passing to extensions $k(f^{\wt \mu_r})$
with truncated Beltrami coefficients, whose values $\vk_p(f^{\wt \mu_r})$ may dominate  $\vk_p(f^\mu)$.

We show here that the equality (13) (i.e., the original version of Grinshpan's conjecture) is valid
for univalent functions $f$ mapping the disk onto domains with asymptotically conformal  boundaries
$L$, which means that for any pair of points $a, b \in L$, we have
$$
\max\limits_{z \in L(a,b)} \frac{|a - z| + |z - b|}{
|a - b|} \to 1 \quad \text{as} \quad |a - b| \to 0,
$$
where the point $z$ lies on $L$ between $a$ and $b$.

Such curves are quasicircles without corners and can be rather pathological
(see, e.g., \cite{Po2}, p. 249). All $C^1$-smooth curves are asymptotically conformal.

This generalized conformality yields, in particular, that the Schwarzian $S_f$ satisfies
 \be\label{16}
\lim\limits_{|z| \to 1-} (1 - |z|^2)^2 S_f (z) = 0.
\end{equation}

\bigskip\noindent
{\bf Theorem 5}. {\it The equality (13) is valid for any asymptotically conformal function
$f \in \Sigma$ with $\|S_f\|_\B < 2$.
}

\bigskip\noindent
{\bf Proof}. The equality (16) implies the existence of Strebel's frame map \cite{St1} for $f$,
and therefore the function $f$ has (extremal) Teichm\"{u}ller extension $f^{\mu_0}$ to $\D$.
By Ahlfors-Weill \cite{AW}, this function admits also extension with harmonic Beltrami coefficient
$$
\mu_{AW}(z) = - \fc{1}{2} (1 - |z|^2)^2 \ S_f(1/\ov z) (1/\ov z^4), \quad |z| <  1,
$$
and in view of the characteristic properties of Teichm\"{u}ller coefficients (see e.g. \cite{GL}, \cite{Kr1}), we have for such coefficients the equality
 \be\label{17}
\mu_{AW} = \mu_0 + \sigma, \quad  \sigma \in A_1(\D)^\bot,
\end{equation}
where
$$
A_1(\D)^\bot = \{\nu \in \Belt(\D)_1: \ \langle \nu,\psi\rangle_\D = 0 \
\text{for all} \ \psi \in A_1(\D)\} = \ker \phi_\T^\prime(\mathbf 0)
$$
is the set of so-called locally (infinitesimally) trivial Beltrami coefficients.

Now observe that since the quantity $\a(f)$ and the representation (10) in Theorem 1 involve only
the integrals $\iint_\D \mu_0 \psi dx dy$ over the unit sphere in the space $A_1(\D)$, the relation
(17) allows us to replace in these integrals the factor $\mu_0$ by the corresponding $\mu_{AW}$.

In addition, the property (16) of Schwarzians of asymptotically conformal functions provides that
the truncated harmonic coefficients $[\mu_{AW}]_\rho$ satisfy
 \be\label{18}
\|[\mu_{AW}]_\rho - \mu_{AW} \|_\iy \to 0 \quad \text{as} \ \ \rho \to 1
\end{equation}
(and accordingly, $\|S_{f^{[\mu_{AW}]_\rho}} - S_{f^{\mu_{AW}}}\|_\B \to 0$).

All this allows one to repeat the arguments from the above proof of Theorem 4 choosing now the annuli
$\{1 < |z| < 1/\rho_j\}$ in accordance with (18). This leads in similar fashion to relations of type (14).

Now, using the convergence of $S_{f^{[\mu_{AW}]_\rho}}$ to $S_{f^{\mu_{AW}}}$ in $\B$-norm and the domination
of dilatation $k(f)$ over all $\vk(f_p)$, one derives from (14) the desired relation (13), completing
the proof.

\bigskip
Note that we have for any $\mu \in \Belt(\D)_1$ the estimate  $\|S_{f^\mu} \|_\B \le 6 \|\mu\|_\iy$
(see, e.g., \cite{Kr2}, Ch 5), and for extremal $\mu_0$ also the lower bound
$\|S_{f^{\mu_0}} \|_\B \ge \|\mu_0\|_\iy$.

\bigskip\bigskip
\centerline{\bf 4. APPLICATIONS OF THEOREM 4}

\bigskip
The applications presented in this section  illustrate the importance of the limit Grunsky norm giving various sharp bounds for univalent functions with quasiconformal extensions, for basic curvelinear functionals and for the Carath\'{e}odory metric.

\bigskip\noindent
{\bf 4. 1. Ahlfors' problem}.
In 1963, Ahlfors posed in \cite{Ah2} (and repeated in his book \cite{Ah3}) the following question which gave rise to various investigations of quasiconformal extendibility of univalent functions.

\bigskip\noindent
{\bf Question}. {\it Let $f$ be a conformal map of the disk (or half-plane) onto
a domain with quasiconformal boundary (quasicircle). How can this map be characterized? }

\bigskip
He conjectured that the characterization should be in analytic properties of the logarithmic derivative $\log f^\prime = f^{\prime\prime}/f^\prime$,
and indeed, many results on quasiconformal extensions of holomorphic maps have been established using $f^{\prime\prime}/f^\prime$ and other invariants (see, e.g., the survey \cite{Kr7} and the references there).

This question relates to another still not completely solved problem in geometric complex analysis:

\bigskip
{\it To what extent does the Riemann mapping function $f$ of a Jordan domain $D \subset \hC$ determine the geometric and conformal invariants (characteristics) of the boundary $\partial D$ and of complementary domain $D^* = \hC \setminus \ov D$?}

\bigskip
Theorem 1 implies a natural qualitative answer to all these questions and shows how the inner features of
conformality completely prescribe the admissible distortion under quasiconformal extensions
of function $f$ and determine the hyperbolic features of the universal Teichm\"{u}ller space.
We present this important consequence from Theorem 1 as

\bigskip\noindent
{\bf Theorem 6}. {\it For any function $f$ mapping conformally the unit disk onto a domain with quasiconformal boundary $L = f(|z| = 1)$, the reflection coefficient $q_L$ of the curve $L$ is determined by the limit Grunsky norm of this function via
 \be\label{19}
 \fc{1 + q_L}{1 - q_L} = \left(\fc{1 + \wh \vk(f)}{1 - \wh \vk(f)}\right)^2.
\end{equation}
Hence the right-hand side of (19) represents the minimal dilatation of quasiconformal reflections across $L$. }

\bigskip\noindent
{\bf 4.2. Invariant metrics on universal Teichm\"{u}ller space}.
As was mentioned above, the equalities (13) and (15) are naturally connected with the invariant distances on the space $\T$. Here we give their explicit representation generated by the original univalent
functions.

It is elementary that the Carath\'{e}odory, Kobayashi and Teichm\"{u}ller metrics of any Teichm\"{u}ller space $\wt T$ are related by
  \be\label{20}
c_{\wt \T}(\cdot, \cdot) \le d_{\wt \T}(\cdot, \cdot) \le \tau_{\wt\T}(\cdot, \cdot),
\end{equation}
and similarly for their infinitesimal forms. The Royden-Gardiner theorem yields that $d_{\wt \T}$ and $\tau_{\wt \T}$ coincide.

As was mentioned, Theorem 3 intrinsically relates to invariant metrics on the space $\T$ giving explicit  representation of these metrics by $\wh \vk(f)$.

\bigskip\noindent
{\bf Theorem 7}. {\it The Carath\'{e}odory metric of the universal Teichm\"{u}ller space $\T$ coincides with its Teichm\"{u}ller metric; hence all non-expansive holomorphically invariant metrics on the space $\T$ are equal, in particular, for any two point $\vp_1, \ \vp_2 \in \T$,
In particular, for any pair $(\vp_1, \ \vp_2) \in \T \times \T$,
 \be\label{21}
c_\T(\vp_1, \vp_2) = d_\T(\vp_1, \vp_2) = \tau_\T(\vp_1, \vp_2),
\end{equation}
and similarly for the infinitesimal forms of these metrics.   }

\bigskip\noindent
{\bf Proof}. It follows from Theorem 1 and the relations (7), (15), (20) that for any point $S_{f^\mu} \in \T$ there exists a sequence of holomotphic functions
$$
 h_{\x,\rho_j,p_j}(S_f) =  \sum\limits_{m,n=1}^\iy \ \sqrt{m n} \
\a_{m n}(f^{(\mathcal R_{p_j}^* \mu_0)_{\rho_j}}) \  \rho_j^{m+n} x_m x_n: \ \T \to \D,
$$
such that
$$
\lim_{j \to \iy} \tanh^{-1} |h_{\x,\rho_j, p_j}(S_f))|  = \tau_\T(\mathbf 0, S_{f^\mu}).
$$
Equivalently,
 \be\label{22}
c_\T(S_{f^\mu}, \mathbf 0) = d_\T(S_{f^\mu}, \mathbf 0) = \tau_\T(S_{f^\mu}, \mathbf 0)
= \tanh^{-1} \wh \vk(f^\mu).
\end{equation}

Now consider two arbitrary points $\vp_1 = S_{f_1}$ and $\vp_2 = S_{f_2}$ in $\T$.
Since the universal Teichm\"{u}ller space is a complex homogeneous domain in $\B$, this general case is reduced to (22) by moving one of these points to the origin $\vp = \mathbf 0$, applying a right translation of the space $\T$ (which is equivalent to replacing the base point of the space $\T$).
Such translations preserve the invariant distances.

The proof of infinitesimal version of relations (21) involves the results lying in the different frameworks (see, e.g., \cite{Di}, \cite{Kr7}) and will not be considered here.

The equality of Carath\'{e}odory metric of universal Teichm\"{u}ller space $\T$ with its Teichm\"{u}ller metric  was established in other way in \cite{Kr13}.

The situation is different in the case of finite dimensional Teichm\"{u}ller spaces. It was established by Gardiner \cite{Ga} and Markovic \cite{Ma}, that all such spaces of dimension greater than $1$ have points with different Carath\'{e}odory and Kobayashi distances; see also \cite{DM}, \cite{Kra}.

\bigskip\noindent
{\bf 4.3. Pluricomplex Green function of universal Teichm\"{u}ller space}.
The equality (21) also determines the potential features of the space $\T$.
We illustrate this by representation of its {\bf pluricomplex Green function}.

Recall that the pluricomplex Green function $g_D(x, y)$ of a domain $D$ in a complex Banach space $X$
with pole $y$ is defined by
$g_D(x, y) = \sup u_y(x) \quad (x, y \in D)$
followed by the upper semicontinuous regularization
$$
v^*(x) = \lim\limits_{\ve\to 0} \sup_{\|x^\prime - x\|<\ve} v(x^\prime).
$$
The supremum  here is taken over all plurisubharmonic functions
$u_y(x): \ D \to [-\iy, 0)$ such that
$u_y(x) = \log \|x - y\|_X + O(1)$
in a neighborhood of the pole $y$.
Here $\|\cdot\|_X$ denotes the norm on $X$, and the remainder term $O(1)$ is bounded from above.
The Green function $g_D(x, y)$ is a maximal plurisubharmonic  function on
$D \setminus \{y\}$ (unless it is identically $- \iy$).

\bigskip\noindent
{\bf Theorem 8.} {\it For every point $\vp = S_f \in \T$, the pluricomplex Green function $g_\T(\mathbf 0, S_f)$ with pole at the origin of $\T$ is given by
$$
g_\T (\mathbf 0, S_f) = \log \wh \vk(f),
$$
and for any pair $(\vp, \psi)$ of points in $\T$, we have
  \be\label{23}
g_\T(\vp, \psi) = \log \tanh d_\T(\vp, \psi) = \log \tanh c_\T(\vp, \psi) = \log k(\vp, \psi),
\end{equation}
where $k(\vp, \psi)$ denotes the extremal dilatation of quasiconformal maps determining the Teichm\"{u}ller distance
between $\vp$ and $\psi$.   }

\bigskip
The first equality (23) is in fact a special case of the general relation
$$
g_D(x, y) = \log \tanh d_D(x, y),
$$
which holds for any hyperbolic Banach domain $D$ whose Kobayashi metric $d_D$ is logarithmically plurisubharmonic  (cf. \cite{Di}, \cite{Kl}, \cite{Kr5}).

Recently the pluripotential theory of Teichm\"{u}ller spaces was created by Miyachi in \cite{Miy1}, \cite{Miy2}.

\bigskip\noindent
{\bf 4.4. Strengthening Pommerenke-Zhuravlev's theorem}.
As was mentioned above, Theorem 4 allows us to establish the minimal possible bound in Theorem 2.

\bigskip\noindent
{\bf Theorem 9}. {\it Every univalent function $f(z) \in \Sigma$ with Grunsky norm $\vk(f) < 1$
admits quasiconformal extensions $f^\mu$ to the whole sphere $\hC$ with dilatations
 \be\label{24}
\|\mu\|_\iy \ge  \limsup \limits_{p \to \iy} \vk_p(f) = \wh \vk(f).
\end{equation}
The bound $\wh \vk(f)$ sharp and best possible for any $f$.
}

\bigskip\noindent
{\bf Proof}. Let $\vk(f) = k < 1$. By Theorem 2, $f$ admits $k_1(k)$-quasiconformal extension to
$\hC$ with $k_1(k) \ge k$, and hence, its dilatation $k(f) \le k_1$. This yields that all
$f_p(z) = f(z^p)^{1/p}$ also have $k_1$-quasiconformal extensions. Thus one can apply Theorem 4,
getting the bound (24).

\bigskip\noindent
{\bf 4.5. Remark}. It would be interesting to find whether Theorem 4 provides Theorem 2. This requires
the uniform estimate $\sup_p \vk_p(f) < 1$.

\newpage
\centerline{\bf 5. OTHER METHODS FOR DETERMINATION OF ANALYTIC AND}
\centerline{\bf GEOMETRIC QUASIINVARIANTS OF CURVELINEAR POLYGONS}

\bigskip\noindent
{\bf 5.1}.
The above theorems bridge quantitative estimating conformal and quasiconformal maps onto quasiconformal domains with the quasireflection coefficients and Fredholm eigenvalues of the boundary curves $L$ of
the image domains.

One of the main still open problems in applications of quasiconformal analysis is  to create the algorithms for explicit determination of these quantities for individual quasicircles or quasiintervals.

In the general case, one had only the rough upper bound for $1/\rho_L$ is given by Ahlfors' inequality \cite{Ah1}
  \be\label{25}
\fc{1}{\rho_L} \le q_L,
\end{equation}
where $q_L$ denotes the minimal dilatation of quasireflections across $L$.

In view of the invariance of all quantities in (25) under the action of the M\"{o}bius group $PSL(2, \C)/\pm \mathbf{1}$, it suffices to consider the quasiconformal homeomorphisms of the sphere carrying
the unit circle $\mathbb S^1$ onto $L$ whose Beltrami coefficients $\mu_f$ have support in the unit disk.  Then $q_L$ is equal to the minimum $k(f)$ of dilatations $\|\mu\|_\iy$ of quasiconformal extensions of the function $f^* = f|\D^*$ into $\D$, and the inequality (25) is equivalent to $\vk(f) \le k(f)$.

The known results in this field with classical origins are presented, for example, in surveys  \cite{Kr8}, \cite{Kr15}, \cite{Ku3}, \cite{Ku5}.

As was mentioned after Theorem 2, all previous estimates of the associated analytic and geometric
curvelinear quasiinvariants have been obtained only for curves related to univalent functions with $\vk(f) = k(f)$. We disclose here some new recent results.

\bigskip\noindent
{\bf 5.2. Sobolev's Beltrami coefficients}.
Assume that a Beltrami coefficient $\mu_0 \in \Belt(D)_1$ is extremal in its class but not of Teichm\"{u}ller type.
A point $z_0 \in \partial D$ is called {\bf substantial} (or essential) for $\mu_0$ if
for any $\ve > 0$ there exists a neighborhood $U_0$ of $z_0$ such that
$$
\sup_{D^*\setminus U_0} |\mu_0(z)| < \|\mu_0\|_\iy - \ve;
$$
so the maximal dilatation $k(w^{\mu_0}) = \|\mu\|_\iy$ is attained on $D$ by approaching
this point.

In addition, there exists a sequence $\{\psi_n\} \subset A_1(D)$ such that $\psi_n(z) \to 0$
locally uniformly on $D$ but $\|\psi_n\| = 1$ for any $n$, and
$$
\lim\limits_{n\to \iy} \iint\limits_D \mu_0(z) \psi_n(z) dx dy = \|\mu_0\|_\iy.
$$
Such sequences are called {\bf degenerated}.

The image of a substantial point is a common point of two quasiconformal arcs, which can be of
spiral type.

Now consider now the univalent functions $f(z) \in S_Q(D^*)$ whose restrictions to the boundary quasicircle $L = \partial D^*$ have substantial points, hence do not have the Teichm\"{u}ller extremal extensions across $L$.

Fix $p > 1$ and consider the subset $\mathcal M_p(D)$ of $\Belt(D)_1$ whose $\mu$  satisfy:

$(i)$ \ $\mu \in L_\iy(D) \bigcap W^{1,p}(D)$,where $W^{1,p}(D)$ is the Sobolev space of functions $\mu$ in $D$ having the first distributional derivatives which belong to $L_p(D)$;

$(ii)$ the value $\|\mu\|_\iy = \esssup_D |\mu(z)|$ is attained by approaching $z$ the boundary of $D$;

$(iii)$ there is a subarc $\gamma \subset \partial D$ depending on $\mu$ such that $\mu(z) \to 0$ as $z$ approaches $\gamma$
from inside $D$.

The boundary values $\mu(z_0)$ for $z_0 \in \partial D$ must be understand as $\lim\limits_{z \to z_0} \mu(z)$.
Note that the value  $\|\mu\|_\iy $ also can be attained by approaching the inner points of domain $D$ and that the arc $\gamma$ is
locally $C^\a$ smooth with $\a > 0$ depending on $\|\mu\|_\iy$, in accordance with the H\"{o}lder continuity of quasiconformal automorphisms of $\hC$.

\bigskip\noindent
{\bf Theorem 10}. \cite{Kr15} {\it For any $p > 2$, every Beltrami coefficient $\mu \in \mathcal M_p(D)$ is extremal in its equivalence class $[\mu]$, and the corresponding quasiconformal automorphism $f^\mu$ of $\hC$ satisfies}
 \be\label{26}
k(f^\mu) = \vk_{D^*}(f^\mu) = \|\mu\|_\iy.
\end{equation}

\bigskip
In the case of $D = \D$, the equalities (26) are completed by terms containing the reflection
coefficient $q_{f(\mathbb S^1)}$ and Fredholm eigenvalue $\rho_{f(\mathbb S^1)}$ similar to Theorem 2.

The situation described by this theorem does not appear in the case of sufficiently
high boundary regularity of univalent functions $f$ (and of $\partial D$) because, for example, for
any  $C^{2+\a}$ smooth $\mu$ the map $f^\mu$ is $C^{2+\a}$ on $\C$), and by Strebel's frame mapping criterion \cite{GL}, \cite{St1}) the equivalence class of $f^\mu$ contains unique Teichm\"{u}ller coefficient $\mu_0 = k |\psi_0|/\psi_0$ with $\psi_0 \in A_1(D)$, which cannot have the substantial boundary points.

The assumptions on the set $\mathcal M_p$ are natural and cannot be replaced in terms of smoothness or non-smoothness of $\mu$ on the boundary points. The same concerns the assumtion $p >2$, which is essential in the proof. It cannot be weakened without adding some additional conditions.

Similar results are valid for harmonic Beltrami coefficients
$$
\mu_\vp(z) = \lambda_D^{-2}(z) \ov{\psi(z)}, \quad \vp \in \B(D),
$$
(in the Kodaira-Spencere sence) defining quasiconformal automorphisms $f$ of the sphere $\hC$ whose
restrictions to $D^*$ are holomorphic univalent with $\vp = S_f$, and such that the maximum of the function $\ld_D^{-2}(z) \ov{S_f(z)} $ is attained at some boundary point and $S_f$ is bounded on some subarc $\gamma$ of $\partial D^*$ .
Here $\ld_D(z) |dz|$ is the hyperbolic metric of domain $D$ with Gaussian curvature $- 4$.

Note that Theorem 10 also provides the unique extremality for some classes of Beltrami coefficients different from Teichm\"{u}ller type.

\bigskip\noindent
{\bf 5.3}.  The equalities (26) are obtained by applying the conformal map $\chi$ of domain $D$ onto the
half-strip
$$
\Pi_{+} = \{\zeta = \xi + i \eta : \ \xi > 0, \ 0 < \eta < 1\}
$$
onto $D$ moving the point $z_0$ of maximality $|\mu(z)|$ to $\iy$ and the endpoints of the
arc $\gamma$ to $0$ and $1$. This pulled the coefficient $\mu$ back to Beltrami coefficient
$$
\nu(\z) := \chi^* \mu(z) = (\mu \circ \chi^{-1})(\z) \ \chi^\prime(\z)/ \ov{\chi^\prime(\z)}
$$
on $\Pi_{+}$. The assumptions on the coefficient $\mu(z)$ and the appropriate smoothness of conformal
map $\chi$ provide the equality
\be\label{27}
\lim\limits_{m\to \iy} \iint\limits_{\Pi_{+}} \nu(\z) \omega_m(\z) d \xi d \eta
= \int\limits_0^1 d \eta \ \lim\limits_{m \to \infty} \fc{1}{m} \int\limits_0^\iy
\nu(\xi + i \eta) e^{- \xi/m} d \xi = \nu(\iy),
\end{equation}
where
$$
\omega_m(\z) = \frac{1}{m} e^{- \z/m}, \quad \z \in \Pi_{+} \ \ (m = 1, 2, ...)
$$
is a degenerating sequence for the affine horizontal stretching of $\Pi_{+}$, which transforms into a
degenerating sequence for the original Beltrami coefficient $\mu$ on $D$.
The relation (27) implies the extremality of $\mu$ and the desired equalities (26) (all  details are
given in \cite{Kr14}).

\bigskip\noindent
{\bf 5.4}. Surprisingly this approach found application to solving an old problem in Teichm\"{u}ller space theory (see \cite{Kr19}).

\bigskip\noindent
{\bf 5.5. Convex curvelinear polygons with infinite number of vertices}. Another approach to explicit calculation of the indicated above curvelinear functionals involves the intrinsic Finsler metrics of negative holomorphic curvature and requires a comparison of such metrics.
We illustrate this on infinitely sided polygons. The case of curvelinear polygons with  finite number of
sides is described in \cite{Kr14}.

Consider an unbounded curvelinear polygon $P_L$ with a countable set of vertices and whose sides are the closures of smooth quasiconformal intervals such that their union $L$ is a quasicircle, and the endpoints of these quasiintervals can accumulate only to a countable set on the boundary $L$ including the infinite point. Let $f$ maps conformally
the disk $\D$ onto $P_L$ and $f^*$ maps $\D^*$ onto the complementary polygon
$P_L^* = \hC \setminus \ov{P_L}$. Denote the vertices of angles of $P_L$ by $A_j = f(z_j), \ z_j \in \mathbb S^1$, their openings by $\pi \a_j, \ 0 < \a_j < 1$, and let the angle in the point $A_\iy = f(z_\iy) = \iy$ be equal $- \pi \a_\iy$ with $\a_\iy > 0$.
Assume that
$$
0 < \a_0 \le \sup_{j \in \mathbb N} |\a_j| <1.
$$
For such polygons, we have the following estimates.

\bigskip\noindent
{\bf Theorem 11}. \cite{Kr17} {\it Under the above assumptions, both conformal mapping functions
$f: \ \D \to P_L$ and $f^*: \ \D^* \to P_L^*$ have equal Teichm\"{u}ller and Grunsky norms, and
 \be\label{28}
\vk(f)  = \vk(f^*) = k(f) = k(f^*) = q_L = 1/\rho_L \ge 1 - |\a|,
\end{equation}
where
$$
1 - |\a| = \max \ \bigl(\sup_{j \in \mathbb N} (1 - \a_j), \ 1 - |\a_\iy|\bigr).
$$
The last inequality in (28) is sharp. In the case $|\a_\iy| < a_j$ for any $j$,  all six relations in (28) are the equalities.

A similar assertion is valid for the unbounded concave domains (not containing inside the infinite point); for those one must replace the lower bound in (28) by $|\beta| - 1$, where $\pi |\beta|$ is equal to supremum of openings of the interior angles of $P_L^*$.  }

\bigskip
The proof of this theorem is based on comparison of the Finsler metric $\ld_\vk$ generated by
the Grunsky coefficients ( generalized Gaussian curvature $\kappa_{\ld_\vk} \le - 4$ with the infinitesimal  form of Teichm\"{u}ller-Kobayashi metric on appropriate holomorphic disk in $\T$.

\bigskip\noindent
{\bf 5.6. Illustrating examples}.
The following examples illustrate Theorem 11 and give rise to a new direction in this field.

\bigskip\noindent
{\bf Example 2} {\it (rectilinear polygon with infinite number of vertices)}.
It follows from Theorem 1 that for any (oriented) closed unbounded curve $L$  with the convex interior which is $C^{1+ \dl}$ smooth at all finite points and has at infinity the asymptotes approaching the interior angle $\pi \a_0 <0$, we have
the equalities
  \be\label{29}
\vk(f) = k(f) = q_L = 1/\rho_L = 1 - |\a_0|,
\end{equation}
where $f$ maps conformally the disk $\D$ onto the inner domain of $L$.

We approximate the curve $L$ by an unbounded rectilinear polygonal line $L^\prime$ going to infinite point in both directions and
such that the inner angles between any two neighbour intervals have the openings $\pi \a_n \in (0, 1)$ satisfying $\a_n > \a_0$ for every $n \in \mathbb N$.

It follows from Theorem 11 that for any such line the values of its curvelinear quasinvariants are explicitly given similar to (29) by
$$
\vk(f_\iy) = k(f_\iy) = \vk(f_\iy^*) = k(f_\iy^*) = q_{L_\iy} = 1/\rho_{L_\iy} = 1 - |\a_0|,
$$
where $f_\iy$ and $f_\iy^*$ are, respectively, the interior and the exterior conformal mapping functions of the
infinite polygonal line $L^\prime$.

\bigskip\noindent
{\bf Example 3} {\it (reflections across subsets of infinite polygonal arcs).} The above example simultaneously implies
the exact estimate for quasiconformal reflections across proper subsets $E$ of polygonal arcs (defined as the orientation
reversing  quasiconformal involutions of the sphere $\hC$ preserving these sets fixed). The reflection coefficients $q_E$
of such sets are determined similar to the case of curves.

Consider a convex unbounded polygon $P_L$ whose boundary $L$ is the union of countably many rectilinear segments (the sides of $P_L$) and the set of vertices $A_j$ of $P_L$ accumulating to the infinite point $A_\iy = \iy$.
The opening $\pi \a_\iy$ at $A_\iy$ is understanding as the infimum of (negative) openings of angles adjacent to $A_\iy$ from inside of $P_L$  and for which the intervals $(A_j, \iy)$ are one of their sides.
Assume that the angles $\pi \a_j$ at all finite vertices $A_j$ satisfy $\a_j < |\a_\iy|$.

In contrast to the case of (closed) curves, only a few estimates have been established for quasireflections across the arcs.
We mention here that the features of quasireflections and convexity of polygon $P_L$ provide the following
interesting phenomenon: take, for example, a finite collection of sides
$$
\mathcal A^{(m)} = \{A_{j_1}, \dots, A_{j_m}\}
$$
and consider the union $L^\prime = L \setminus \mathcal A^{(m)}$; then every set of the form
$$
E = L^\prime \cup e,
$$
where $e$ is a subset of $\mathcal A^{(m)}$, has the reflection coefficient $q_E = 1 - |\a_\iy|$,
which is equal to the reflection coefficient of the whole polygon $P_L$.

\bigskip\noindent
{\bf Example 4} ({\it Rational and harmonic Beltrami coefficents}).
Any rational function $r_n \in \B$ with poles of order two on the unit circle of the form
$$
r_n(z) = \sum\limits_1^n \fc{c_j}{(z - a_j)^2} + \sum\limits_1^n \fc{c_j^\prime}{z - a_j},
$$
satisfying $\sum\limits_1^n |c_j| > 0$ and $\|r_n\|_\B < 2$, generates the harmonic Beltrami
coefficients
$$
t \nu_{r_n} = - \fc{t}{2} (1 - |z|^2)^2 \ r_n(\ov z), \quad |t| < 1,
$$
whose moduli attain their maximal values $|t|$ on the boundary. So, all these coefficients have
substantial points on $\mathbb S^1$, and the above arguments imply the equalities
$$
\vk(f^{t \nu_{r_n}}) = k(f^{t \nu_{r_n}}) = q_{f^{t \nu_{r_n}}} = 1/\rho_{f^{t \nu_{r_n}}}.
$$
The same is also valid for more general harmonic coefficients.
The equivalence classes of all such Beltrami coefficients do not contain the Teichm\"{u}ller
coefficients.

\bigskip\bigskip
\centerline{\bf 6. EXTREMALITY OF GRUNSKY-MILIN FUNCTIONAL}

\bigskip\noindent
{\bf 6.1. General remarks}.
Consider again a bounded quasicircle $L \subset \C$ with the interior $D \ni 0$ and exterior $D^* \ni \iy$ and the corresponding classes $S_D(\iy)$ and $\Sigma_{D^*}(0)$ and investigate on these classes
the bounded holomorphic (continuous and Gateaux $\C$-differentiable) functionals $J(f)$, which means
that for any $f$ and small $t \in \C$,
$$
J(f + t h) = J(f) + t J_f^\prime(h) + O(t^2), \quad t \to 0,
$$
in the topology of uniform convergence on compact sets in $D$ or $D^*$.
Here $J_f^\prime(h)$ is a $\C$-linear functional which is lifted to the strong (Fr\'{e}chet) derivative of $J$ in the norms  of both spaces $L_\iy(D^*)$ and $\B(D)$ (accordingly, in $L_\iy(D)$ and $\B(D^*)$.
For definiteness, consider again the class $\Sigma_{D^*}(0)$ and its subclasses $\Sigma_{D^*,k}$
of functions with $k$-quasiconformal extensions to $\hC$ with $f(0) = 0$. The union
$\Sigma^0(D^*) = \bigcup_{k < 1} \Sigma_{D^*,k}$ is dense in $\Sigma(D^*)$ in the topology of locally inform convergence on $D^*$.

One can assume, without loss of generality, that
$$
M = \max_{\Sigma(D^*)} |J(f)| = 1,
$$
passing if needed to functional $J^0(f) = J(f)/M$
(the existence of extremal functions maximizing $J$ follows from compactness of $\Sigma(D^*)$ with respect to locally uniform convergence).

Varying $f$, one gets the corresponding functional derivative
 \be\label{30}
\psi_0(z) = J_{\id}^\prime(g(\id, z)),
\end{equation}
where
$$
g(w, \z) = \fc{1}{w - \z}
$$
is the kernel of variation on $\Sigma^0(D^*)$. Any such functional is represented by a complex Borel measure on $\C$ and extends thereby to all holomorphic  functions on $D^*$ (cf. \cite{Sch}). We assume that this derivative is meromorphic on $\C$
and has in the domain $D$ only a finite number of the simple poles (hence $\psi_0$ is integrable over $D$).

For example, one can take the general distortion functionals
 \be\label{31}
J(f) = J(f(z_1), f'(z_1), \dots \ , f^{(\a_1)}(z_1); \dots; f(z_p), f'(z_p), \dots \ , f^{(\a_p)}(z_p))
\end{equation}
with $\grad \wh J(\mathbf 0) \ne 0$, where $z_1, \dots \ , z_p$ are distinct fixed points in $D^*$ with assigned orders $\a_1, \dots, \a_p$, and $J$ is a holomorphic function of its arguments.

In this case, $\psi_0$ is a rational function
$$
\wh J_{\id}^\prime(g(\id, z)) = \sum\limits_{j=1}^p \sum\limits_{s=0}^{\a_j-1} \ \fc{\partial
\wh J(\mathbf 0)}{\partial \om_{j,s}} \fc{d^s}{d\z^s} g(w,\z)|_{w=z,\z=z_s},
$$
where $\om_{j,s} = f^{(s)}(z_j)$.

We associate with the derivative (29) the Teichm\"{u}ller disk
$$
\D(\mu_0) = \{t \mu_0: \ |t| < 1, \ \ \mu_0 = |\psi_0|/\psi_0\}.
$$

\bigskip\noindent
{\bf 6.2. General main theorem}.
For all such functionals, we have the following two general theorems, which provide a useful tool for explicit solving the extremal problems for univalent functions and essentially improves the known methods and results.

\bigskip\noindent
{\bf Theorem 12}. {\it If all zeros of the functional derivative $\psi_0(z) = J_{\id}^\prime(g(\id, z))$ on $D$ of a given functional (31) are of even order and $\psi_0(0) \ne 0$, then for any $k, \ 0 < k < 1$,
the maximal value of this functional on the any class $\Sigma_k(0)$ is attained at the point ${k|\psi_0|/\psi_0}$ on the disk $\D(\mu_0)$, and
 \be\label{32}
\max_{\Sigma_k(0)} |J(f^\mu)| = |J(f^{k|\psi_0|/\psi_0})| = k = \a(f),
\end{equation}
with $\a(f)$ given by (11).

In addition, the indicated maximal value of $|J(f^\mu)|$ is attained only when this functional coincides  on the disk $\D(\mu_0)$ with the Grunsky-Milin functional
$$
J_{\x^0}(f^\mu) =  \sum\limits_{m,n=1}^\iy \beta(f^\mu) x_m^0 x_n^0
$$
whose defining point $\x^0 = (x_n^0) \in S(l^2)$ is determined by $\mu_0$.
}

\bigskip
Recall that the Grunsky-Milin coefficients of functions $f \in \Sigma(D^*)$ are determined from expansion (3).

In the case of the canonical dick $\D^*$, this theorem  can be essentially  strengthened, getting
simultaneously the values of the basic analytic and geometric quasiinvariants of the extremal curves
$L_0 = f^{\mu_0}(\mathbb S^1)$.

\bigskip\noindent
{\bf Theorem 13}. {\it Let a functional $J(f)$ be considered on the class $\Sigma_Q$ and its
functional derivative $\psi_0(z) = J_{\id}^\prime(g(\id, z))$  has in the unit disk only zeros of even order. Suppose also that $\psi_0(0) \ne 0$.
Then the maximal value of the modulus of $J(f^\mu)$ on any subclass $\Sigma_k(0), 0 <k < 1$, is attained at $\mu_0(z) = {k|\psi_0|/\psi_0}$.

In addition, the reflection coefficient and the Fredholm eigenvalue of the curve
$L_1 = f^{\mu_0}(\mathbb S^1)$ are explicitly given by
 \be\label{33}
q_{L_1} = 1/\rho_{L_1} = \a(f).
\end{equation}
}

\noindent
{\bf 6.3.} The proof of Theorem 12 also involves Theorem 1.
Both functionals $J(f^\mu)$ and the Grunsky functional have the same differential at the origin of the disk $\D(\mu_0)$ and, since either from these functionals maps the unit disk into itself, the Schwarz lemma yields that these functionals must coincide on the entire disk $\D(\mu_0)$.

Theorem 12 implies (since $\psi_0(0) = \grad J_{\id} (\mathbf 0) \ne 0$) that any value attained by  $|J(f^{t \mu_0/\|\mu_0\|_\iy})|$ on the disk $\D(\mu_0$) in $\Belt(\D)_1$ is extremal for all
$|t| \le t_0$ with  fixed $t_0 < 1$.

To compare these value with values of $J(f^\mu)$ on other holomorphic disks in $\Belt(\D)_1$, note that any holomorphic functional $J(f)$ generates a holomorphic map $\eta_1$ of the ball $\Belt(\D)_1$ into
the unit disk via $\eta(\mu) = J(f^\mu)$,
and the Schwarz lemma yields for any $\mu \in \Belt(D)_1$ the inequality
\footnote{This inequality is given also by Lehto's majoration principle}
$$
|J(f^\mu)| \le \|\mu\|_\iy,
$$
which implies the assertion of Theorem 12 on maximality of $|J(f^{\mu_0})|$.

It remains to establish that for all $|t| <1$ these maximal values $J(f^{t\mu_0})$ are equal to values of the Grunsky-Milin functional  $J_{\x^0}(f^{t\mu_0})$. This is obtained by applying the variational formula
$$
f^\mu(z) = z - \fc{1}{\pi} \iint\limits_D \fc{\mu(w)}{w - z} du dv + O(\|\mu\|_\iy^2)
$$
with small $\|\mu\|_\iy$. Its kernel is represented for $z \in D^*$ in the form
$$
\fc{1}{w - z} = \sum\limits_1^\iy P_n^\prime(w) \vp_n(z),
$$
where $\vp_n = \chi^n$ are given in (4) and $P_n$ are well-defined polynomials; the degree of $P_n$
equals $n$.
These polynomials satisfy
$$
\langle P_m^\prime, P_n^\prime \rangle_D = \pi \dl_{m n},
$$
which means that the polynomials $P_n^\prime(z)/\sqrt{\pi}$ form an orthonormal system
in $A_1^2(D)$, and this system is complete (the details see in \cite{Kr10}, \cite{Mi}).
Hence,
$$
\fc{f^\mu(z) - f^\mu(\z)}{z - \z} = 1 - \fc{1}{\pi} \iint\limits_D \fc{\mu(w) du dv}{(w - z)(w - \z)} + O(\|\mu^2\|_\iy)
$$
and (3) implies
$$
\begin{aligned}
- \log \fc{f^\mu(z) - f^\mu(\z)}{z - \z}
&= - \log \Bigl[ 1 - \fc{1}{\pi} \iint\limits_D
\fc{\mu(w) du dv}{(w - z)(w - \z)}\Bigr] + O(\|\mu^2\|_\iy) \\
&= \fc{1}{\pi} \iint\limits_D
\fc{\mu(w) du dv}{(w - z)(w - \z)} + O(\|\mu^2\|_\iy) \\
&= \fc{1}{\pi} \iint\limits_D \mu(w)
\sum\limits_1^\iy P_m^\prime(w) \vp_m(z)
\sum\limits_1^\iy P_n^\prime(w) \vp(\z) du dv + O(\|\mu^2\|_\iy),
\end{aligned}
$$
where the ratio $O(\|\mu^2\|_\iy)/\|\mu^2\|_\iy$ is uniformly bounded on compact sets of $\C$.
Comparison with the representation
$$
- \log \fc{f^\mu(z) - f^\mu(\z)}{z - \z}
= \sum\limits_1^\iy \beta_{m n} \vp_m(z) \vp_n(\z) \quad (\vp_n = \chi^n)
$$
yields
$$
\beta_{m n}(f^\mu) = - \fc{1}{\pi} \iint\limits_D \mu(w) P_m^\prime(w) P_n^\prime(w) du dv
+ O(\|\mu^2\|_\iy),
$$
which provides the representation of differentials of holomorphic functions
$\mu \mapsto \wh \beta_{m n}(\mu) = \beta_{m n}(f^\mu)$ on $\Belt(D)_1$ at the
origin. Using the estimate (5) ensuring the holomorphy of the corresponding functions (4) on this ball,
one obtains that the differential of $h_\x(\mu_0)$ at zero is represented in the form
$$
d h_\x(\mathbf 0) \mu = - \fc{1}{\pi} \iint\limits_D \mu_0(z) \sum\limits_{m,n=1}^\iy
x_m x_n \ P_m^\prime(z) P_n^\prime(z) dx dy, \quad \x = (x_n) \in S(l^2).
$$
This corresponds the case of equality in (10).

The differential of $J(f^{t\mu_0}$ at the origin has the same expression, and the equality
$$
J(f^{t\mu_0}) = J_{\x^0}(f^{t\mu_0}) =  \sum\limits_{m,n=1}^\iy \beta(f^{t\mu_0}) \x_m^0 x_n^0
$$
follows from the Schwarz lemma. This completes the proof of Theorem 12.

As for Theorem 13, its first statement and the right-hand equality in (33) follow from Theorem 12 and from the K\"{u}nau-Schiffer theorem.
Finally, due to \cite{Ah2}, \cite{Ku1} (see also \cite{Kr7}), for any quasicircle $L \subset \hC$, its quasireflection coefficient $q_L$ is connected with the minimal dilatation $k_L = \inf \|\partial_{\ov z} f_{*}/\partial_z f_{*}\|_\iy$
among all orientation preserving quasiconformal automorphisms
$f_{*}$ of $\hC$ carrying the unit circle onto $L$ by
$$
\fc{1 + q_L}{1-q_L} = \Bigl(\fc{1 + k_L}{1- k_L}\Bigr)^2.
$$
This proves the first equality in (33).

Note that Theorems 12 and 13 do not include the uniqueness of extremal functions.

Similar theorems are also valid for univalent functions on the interior (bounded) domain $D$ with expansions  $f(z) = z + a_2 z^2 + \dots$ near the origin.

\bigskip\noindent
{\bf 6.4. Two examples}. Theorem 12 and its counterpart for $S_Q(D)$ provide various sharp distortion estimates
for biholomorphic maps of {\bf arbitrary} quasiconformal domains. In particular, let
$\om(z) = p(z)/q(z)$ be a rational function whose poles (zeros of $q(z)$) are simple and located
on a
bounded quasicircle $L$ with the interior domain $D \ni 0$ and exterior domain $D^* \ni \iy$, and
let $p^\prime(z) \ne 0$ on $D$ (more generally, $p$ has in $D$ only zeros of even order).
Letting
$$
\mu_{*}(z) = |\om(z)|/\om(z), \quad z \in D,
$$
and extending $\mu_{*}$ by zero to $D^*$, consider the solutions $w_k(z)$ of the Beltrami equations
$$
\partial_{\ov z} w = k \mu_{*}(z) \ \partial_z w, \quad 0 < k < 1,
$$
on $\C$ with $w_k(0) = 0, \ w_k^\prime(0) = 1, \ w_k(\iy) = \iy$.
Then by Theorem 1, $\vk_{D^*}(w_k) = k(w_k) = k$ for all such $w_k(z)$.

If $L = \mathbb S^1$, the quantity $k$ is equal to the reflection coefficient of for all $|t| <1$
and is reciprocal to the first positive Fredholm eigenvalue of this curves.

In similar fashion, every polynomial $p(z) = \sum\limits_0^m a_n z^n$ with real zeros, being
composed with a conformal map of the upper half-plane onto $\D$, generates a continuous collection
of quasicircles with quantitatively estimated basic analytic and geometric characteristics.

Another collection of examples of maps with given quasiinvariants is simply obtained
from Theorem 12 as follows.

By this theorem, any $\psi \in A_1^2(D)$ determines for any $k < 1$ the univalent functions
$f^{k|\psi|/\psi} \in \Sigma_Q(D^*)$ with $\vk_{D^*}(f^{k|\pi|/\psi}) = k$.

All this holds for $\psi = \phi^2$, where $\phi$ is a holomorphic function on $D$
with $\phi(z) = o(1/\sqrt{1-|z|})$ as $|z| \to 1$, and in particular, for bounded holomorphic
$\phi$ in $D$.

In the case $D = \D$, one obtains additionally from Theorem 3 the reflection coefficient and Fredholm
eigenvalue of the images of unit circle under these maps.

\bigskip\noindent
{\bf 6.5 Approximate approach}.
In the case of small dilatations, one can use the variational approach.
Let the function $\phi$ be $C^1$-smooth on the closed domain $\ov D$ and let the boundary
$L = \partial D$ be $C^{1+\sigma}$-smooth ($\sigma > 0$). Then the Belinskii theorem
\cite{Bel} implies that for sufficiently small $|t| > 0$, the map
$$
w_\phi(z;t) = z - \fc{t}{\pi} \iint\limits_D \fc{\ov{\phi(\z)}/\phi(\z)}{\z - z} d \xi d \eta
$$
provides a quasiconformal automorphism of the sphere $\hC$ with the complex dilatation
$$
\mu(z, t) = t \ov{\phi(z)}/\phi(z) + O(t^2) \quad \text{for} \ \ z \in D
$$
and conformal on $D^*$, and (10) yields that its dilatation $k(w_\phi(\cdot;t))$ and the Grunsky norm
satisfy
$$
k(w_\phi(\cdot;t)) = \vk_{D^*}(w_\vp(\cdot;t)) = |t| + O(|t|^2),
$$
with uniform estimate of the remainder.

 In the case of arbitrary dilatations and when the functional derivative $\psi_0(z) =  J_{\id}^\prime(g(\id, z))$ has also zeros of odd order in $D$, Theorem 11 can be applied to approximate
solving the extremal problem.
We consider here the case when the number of such zeros is finite.

We shall use the following lemma proven in \cite{Kr1}.
Consider in the space $L_p(\C), \ p > 2$, the well-known integral operators
$$
T\rho(z) = - \fc{1}{\pi} \iint\limits_\C \fc{\rho(\z) d \xi d\eta}{\z - z}, \quad
\Pi \rho(\z) = - \fc{1}{\pi} \iint\limits_\C \fc{\rho(\z) d \xi d\eta}{(\z - z)^2}
= \partial_z T\rho(z)
$$
assuming for that $\rho$ has a compact support in $\C$. Then the second integral
exists as a Cauchy principal value, and the derivative $\partial_z T$ generically is understanding as distributional. One of the basic fact in the theory of quasiconformal maps
is that every quasiconformal automorphism $w^\mu$ with $\|\mu\|_\iy = k < 1$ of the extended plane
$\hC$  with $\|\mu\|_\iy = k < 1$ is represented in the form
$$
w^\mu(z) = z + T\rho(z),
$$
where $\rho$ is the solution in $L_p$ (for $2 < p < p_0(k)$)  of the integral equation
$$
\rho = \mu + \mu \Pi \rho,
$$
given by the series
$$
\rho = \mu + \mu \Pi \mu + \mu \Pi \mu(\mu \Pi(\mu)) + \dots.
$$
Let $B_{p,R}$ denote the space of functions $f(z)$ with $f(0) = 0$ on the disk $\D_R = \{|z| < R\}$ with norm
$$
\|f\|_{B_{p,R}} = \sup_{z_1,z_2 \in \ov{\D_R}}  \ \fc{|f(z_1) - f(z_2)|}{|z_1 - z_2|^{1-2/p}}
+ \|\partial_z f\|_{L_p} + \|\partial_{\ov z} f\|_{L_p}.
$$

\bigskip\noindent
{\bf Lemma 2}. \cite{Kr1} {\it Let $f^\mu$ be a quasiconformal automorphism of $\hC$ conformal in the
disk $\D_R^* = \{|z| > R\}$ normalized by $f^\mu(z) = z + \const + O(1/z)$ as $z \to \iy$ and $f^\mu(0) = 0$. Suppose that $\mu$ satisfies
$$
\|\mu\|_\iy = k < 1, \quad \|\mu\|_{L_r} < \ve,
$$
where
$$
r \ge r_0(k, p) = p_0/ p(p_0 - p),
$$
with $p, p_0 > 2$ indicated above, and $\ve$ is small. Then $f^\mu$ is represented in the form
$$
f^\mu(z) = z - \fc{1}{\pi} \iint\limits_{\D_R} \mu(\z) \Bigl(\fc{1} {\z - z} - \fc{1}{\z}\Bigr)
d \xi d \eta + \om(z, \mu),
$$
where
$$
\|\om\|_{B_p(\ov{\B_R})} \le M(k, p, R) \ve^2,
$$
and the constant $M(k, p, R)$ depends only on $k, p, R$. }

\bigskip
In other words, the bounded Beltrami coefficients, which are integrally small, define quasiconformal variations of the standard form, with uniform estimate of the remainder terms.

Now, let $z_1, \dots, z_s$ be all zeros of $\psi_0(z)$ of odd order in $D \setminus \{0\}$ and $\psi_0(0) \ne 0$.
Fix a small $\ve > 0$ and take a quasidisk $D_\ve \subset D$ with $\mes_2 (D \setminus D_\ve) < \ve$,
where $mes_2$ denotes the $2$-dimensional Lebesque measure on $\C$.

Lemma 2 (reformulated for our maps $f^\mu$ with $\mu$ supported on a bounded domain $D$) yields for any given functional $J(f)$ the uniform estimate
 \be\label{34}
\max_D |J(f^\mu)| - \max_{D_\ve} |J(f^\mu)| = O(\ve) = \const \ \ve,
\end{equation}
and the term $\max_{D_\ve} |J(f^\mu)| = |J(f^{\mu_0})|$ is explicitly described by Theorem 12, which provides the values of $\vk_{D_\ve}(f^{t \mu_0}) |t|$ for all $|t| < 1/\|\mu_0\|_\iy $.

\bigskip\noindent
{\bf 6.6. Examples}. Let $D$ be the {\bf ellipse} $\mathcal E$ considered in \textbf{2.5}, and assume that the zeros of
$\psi_0$ of odd order are placed near the boundary of $\mathcal E$ so that the approximating domain
$\mathcal E_\ve$ also is an ellipse with $\ve$-close semiaxes, where
 \be\label{35}
\ve = \max_j \inf_{\z \in L} |z_j - \z|.
\end{equation}
Then the extremal Beltrami coefficient $\mu_0$ for $\mathcal E_\ve$ is represented by the Chebyshev polynomials from Example 1 for this ellipse in the form
$$
\mu_0(z) = k \ov{\sum\limits_0^\iy x_n^0 P_n(z)} \Big/ \sum\limits_0^\iy x_n^0 P_n(z)
$$
with some $\x^0 = (x_n^0) \in S(l^2)$.

This approximation by ellipses yields a rough explicit accuracy estimate. One can approximate the punctured ellipse $\mathcal E \setminus \{z_1, \dots, z_s\}$ by (unknown explicitly) quasidisks $D_{\ve^\prime}$ satisfying
$\mes_2(\mathcal E \setminus D_{\ve^\prime}) < \ve^\prime$ for arbitrarily small prescribed $\ve^\prime$.

\bigskip
Similar situation is in the case of the {\bf Bernoulli lemniscate}, which is defined by the equation $|z^2 - 1| = 1$ and bounds two domains symmetric with respect to the imaginary axes $i \R$.
The function $w = z^2 - 1$ maps conformally either of these domains onto the unit disk; hence, an
orthonormal basis in $A_2$ is formed by polynomials
$$
P_n(z) = 2 \sqrt{(n + 1)/\pi} \ z (z^2 - 1)^n, \quad n = 0, 1, 2, \dots.
$$

\bigskip
The estimate of type (34) is also valid for more general domains, for example, for starlike quasidisks
$D \subset \C$. An orthonormal basis $\{\psi_n\}$ in $A_2(D)$ is
$$
\psi_n(z) = \sqrt{(n + 1)/\pi} \ h(z)^n h^\prime(z), \quad n = 0, 1, 2, \dots,
$$
where $h$ is a conformal map $D \to \D$. In this case, the approximating domain $D_\ve$ is obtained from
$D$ by stretching $z \mapsto \ve z, \ z \in D$ with $\ve$ given by (35).

\bigskip\bigskip
\centerline{\bf 7. MODELING UNIVERSAL TEICHM\"{U}LLER SPACE BY GRUNSKI}
\centerline{\bf COEFFICIENTS}

\bigskip
There are several models of the universal Teichm\"{u}ller space. The most applicable is the Bers
model via the domain $\T$ in the Banach space $\B$ of the Schwarzian derivatives; this model was
applied above. Other models, also very useful, have been constructed, for example, in \cite{B}, \cite{TT}, \cite{Zh2}.

We mention here another model constructed in \cite{Kr13} by applying the Grunsky coefficients of univalent functions in the disk.
In this model, the space $\T$ is represented by a bounded domain in a subspace of $l_\infty$
determined by the Grunsky coefficients.
This domain is biholomorphically equivalent to the Bers domain $\T$.

\bigskip
Consider again the collection $\Sigma$ of univalent nonvanishing functions in the disk $\D^*$ with hydrodynamical normalization
$$
f(z) = z + b_1 z^{-1} + \dots: \D^* \to \hC \setminus \{0\}.
$$
and their Grunsky \ coefficients $c_{m n}$ defined from the expansion
$$
\log \fc{f(z) - f(\z)}{z - \z} = - \sum\limits_{m, n = 1}^\infty
\a_{m n} z^{-m} \z^{-n}, \quad (z, \z) \in (\D^*)^2;
$$
these coefficients satisfy (2).
Note that $\a_{m n} = \a_{n m}$ for all $m, n \ge 1$ and $\a_{m 1} = b_m$ for any $m \ge 1$.

Let $\Sigma_k$ be the subset of $\Sigma$ formed by the functions with $k$-quasiconformal extensions to $\D$, and $\Sigma^0 = \bigcup_k \Sigma_k$.

These coefficients span a $\C$-linear space $\mathcal L^0$ of sequences
$\mathbf a = (\a_{m n})$ which satisfy the symmetry relation
$\a_{m n} = \a_{n m}$ and
$$
|\a_{m n}| \le C(\mathbf c)/\sqrt{m n}, \quad
C(\mathbf a) = \const < \infty \quad \text{for all} \ \ m, n \ge 1,
$$
with finite norm
 \be\label{36}
\|\mathbf a\| = \sup_{m,n} \ \sqrt{m n} \ |\a_{m n}| +  \sup_{\mathbf x = (x_n) \in S(l^2)} \
\Big\vert \sum\limits_{m,n=1}^\iy \ \sqrt{mn} \ \a_{mn} x_m x_n \Big\vert.
\end{equation}
Denote the closure of span $\mathcal L^0$ by $\mathcal L$ and note that the limits of convergent sequences
$$\{\mathbf a^{(p)} = (\a_{m n}^{(p)})\} \subset \mathcal L
$$
in the norm (36) also generate the double series
$$
\sum\limits_{m,n=1}^\infty \ \a_{m n} z^{-m} \zeta^{-n}
$$
convergent absolutely in the bidisk $\{(z, \z) \in \hC^2: \ |z| > 1, \ |\z| > 1\}$.

It is established in \cite{Kr13} that {\it the sequences $\mathbf a$  corresponding to functions
$f \in \Sigma^0$ with quasiconformal extensions to $\hC$ fill a bounded domain
$\wt \T$ in the indicated Banach space $\mathcal L$ containing the origin, and the correspondence
$$
S_f \leftrightarrow \mathbf a = (\a_{m n})
$$
determines a biholomorphism of this domain $\wt \T$ onto the space $\T$.  }

In this model, the Grunsky norm of any function from $S_Q(\iy)$ arises as a canonical part of the Banach norm of its representative $\mathbf c$ in $\wt \T$, and the hyperbolic length of the limit Grunsky norm
of this function is equal by Theorem 8 to each of invariant distances in $\wt \T$ between the point
$\mathbf c$ and the origin. This determines the basic features of both Grunsky norms.

The corresponding holomorphic functions
$$
 h_{\x^0}(\mathbf a) = \sum\limits_{m,n=1}^\iy \ \sqrt{m n} \ \a_{m n} \ x_m^0 x_n^0
$$
generating the norm $\vk(f)$ become linear on $\wt \T$, which provide some interesting applications.

The basic geometric properties of domain $\wt \T$ are the same as of other models of the space $\T$; in particular, this domain is not starlike, but holomorphically convex, admits pluricomplex potential
description, etc. (cf., e.g., \cite{BE},  \cite{Kr3}, \cite{Kr12}, \cite{Miy1}, \cite{Miy2}).

\bigskip\noindent
{\bf 6.7. Remark}. The maximality of Grunsky-Milin functional gives rise to estimating this functional
on classes $\Sigma(D^*, \vk)$ of univalent functions $f \in \Sigma(D^*)$ with $\vk_{D^*}(f) \le \vk$.
The result is similar to Theorem 12.

\bigskip\bigskip
\centerline{\bf 8. POLYNOMIAL APPROXIMATION OF UNIVALENT FUNCTIONS}
\centerline{\bf WITH QUASICONFORMAL EXTENSION}

\bigskip\noindent
{\bf 8.1. Approximation by univalent polynomials}.
There are deep results in complex geometric function theory that the Riemann conformal mapping function of the disk onto any simply connected hyperbolic domain can be approximated by {\it univalent} polynomials, and generically this approximation is uniform on compact subsets of this disk.

The methods presented in the previous sections also can be applied to get some essential strengthening these results, namely, to obtain the approximation theorems intrinsically connected with quasiconformal extendibility.

First of all, we have for any such function the following result.

\bigskip\noindent
{\bf Theorem 14}. \cite{Kr14} {\it For any univalent function $f$ in the disk $\D$ admitting
quasiconformal extension across the unit circle, there exists a sequence of univalent polynomials $p_n$ on the close disk $\ov \D$ convergent to $f$ uniformly on $\ov \D$ and such that their
dilatations $k(p_n) \nearrow k(f)$.     }

\bigskip
Much stronger results giving approximation in $\B$-norm are valid for asymptotically conformal univalent functions: {\it every univalent function $f(z)$ on $\D$ with assymptotically conformal boundary values  is approximated by univalent polynomials $p_n$ on $\ov \D$ so that $\|S_{p_n} - S_f\|_\B \to 0$ as $n \to \iy$ and $p_n$ admit $k_n$-quasiconformal extensions to $\hC$ with dilatations} $k_n \to k$.

\bigskip\noindent
{\bf 8.2. Quasiinvariants of polynomials}.
Consider the normalized polynomials $p(z) = a_0 + a_1 z + a_2 z^2 + \dots + a_n z^n$
with $a_1 \ne 0$. Their Schwarzians $S_p$ are rational functions on $\hC$ of the form
 \be\label{37}
S_p(z) = \sum\limits_j \fc{c_j}{(z - z_j)^2} \quad \text{with} \ \ \sum\limits_j |c_j| > 0,
\end{equation}
where $z_j$ are the finite critical points of $p$ (i.e., the zeros of its derivative $p^\prime(z)$, which is a polynomial of degree $n - 1$). In the case of univalent $p$ on $\D$, these points are placed outside of $\D$. For such polynomials, we have

\bigskip\noindent
{\bf Theorem 15}. {\it Let a polynomial $p(z)$ be univalent on the disk $\D$ and all finite critical points $z_1, \dots, z_{n-1}$ are placed  on the boundary circle $\mathbb S^1$, and let the $\|S_p\|_\B < 2$. Then
$$
k(p) = \vk(p) = q_{p(\mathbb S^1)} = 1/\rho_{p(\mathbb S^1)} = \max_j |c_j|,
$$
where $c_j$ arise from (37).
}

\bigskip
The proof of this theorem and some its consequences are given in \cite{Kr20}.

All generic features and obstructions arising by the general quasiconformal maps are also
valid for polynomial maps. In particular, there exist univalent polynomials $p_n$ on the closed disk
$\ov \D$ with different Grunsky and Teichm\"{u}ller norms (so, $\vk(p_n) < k(p_n)$).

\bigskip\bigskip
\centerline{\bf 9. GRUNSKY NORM AND BOUNDARY QUASICONFORMAL DILATATION}

\bigskip\noindent
{\bf 9.1}. Consider now the question:

\bigskip
{\it Which intrinsic quasiconformal features of a quasicircle $L = f(\mathbb S^1)$ provide  the equality of its  reflection coefficient to the Grunsky norm of $f$ ? }

\bigskip
The answer will be given not for all quasiconformal curves. It involves the notion of polygonal quasiconformal maps.

Any quasisymmetric homeomorphism $h$ of $\hR = \R \cup \{\infty\}$ onto itself, i.e., such that
$$
M^{-1} \le \frac{h(x + t) - h(x)}{h(x) - h(x + t)} \le M
$$
for any $x \in \R$ and $t > 0$, with $M = M(h) < \iy$, admits quasiconformal
extensions onto the upper half-plane $U = \{z : \Im z > 0\}$
(and symmetrical extensions to the lower half-plane $U^*  = \{z : \Im z < 0\}$).
Let $f_0$ be an extremal quasiconformal extension of $h$.

Now consider quasiconformal maps $f : \ D \to \hC$ of a simply connected
Jordan domain $D \subset \hC$ and regard this domain as a (topological) polygon $D(z_1, ... \ , z_n)$ whose vertices are $n$ distinguished boundary points $z_1, ... \ , z_n$ ordered in accordance with
the orientation of $\partial D$. Assuming that $f$ has a homeomorphic extension to closure $\overline{D}$, one obtains another polygon $D(f(z_1), f(z_2), ... \ , f(z_n))$, and $f$ moves the vertices into vertices.
We shall call such maps {\bf polygonal quasiconformal}.

In view of conformal invariance of $q$-quasiconformality, it will suffice for us to deal with the maps
of polygons whose domains are either the upper half-plane $U$ or the unit disk $\D$.
 maps $f : \ D \to \hC$ of a simply connected Jordan domain $D \subset \hC$ and regard this domain as
 a (topological) polygon $D(z_1, ... \ , z_n)$ whose vertices are $n$ distinguished
boundary points $z_1, ... \ , z_n$ ordered in accordance with the orientation of $\partial D$.

For any quasymmetric map $h$, we have
$k(f_0) = \lim\limits_{n \to \infty} k(f_{[n]})$,
where $f_{[n]}$ are extremal polygonal quasiconformal maps of $U(x_1, ... \ , x_n)$ onto
$U(h(x_1), ... \ , h(x_n))$, and the set of vertices \{$x_n\}$ becomes dense on $\R$.

A long time ago, the question was posed whether for a fixed $n \ge 4$ (or even for $n \le N$) the equality
 \be\label{38}
k(h) = \sup k(f_{[n]})
\end{equation}
is achieved when one allows the vertices to vary on $\hR$ in all possible ways.

In the case $n = 4$, when the polygons are quadrilaterals, the distortion
of conformal modules $\mod Q$ of all quadrilaterals $Q \subset \D$ under a quasiconformal homeomorphism
$f : \ D \to \hC$ determines its maximal dilatation $K(f)$, and it was
conjectured, accordingly, that in this case the best possible bound in (1)
is equal to
 \be\label{39}
K_0(h) = \sup\limits_Q \frac{\mod f_0(Q)}{\mod Q},
\end{equation}
where $f_0$ is an extremal extension of $h$ and supremum is taken over
all quadrilaterals $Q = U(z_1, z_2, z_3 , z_4)$ (i.e., that it suffices
to take only quadrilaterals with the vertices on $\hR$).

This conjecture was disproved by a counterexample of Anderson and Hinkkanen (see \cite{AH}) and by Reich \cite{Re}; a complete answer for quadrilaterals was given in \cite{Wu} (see also \cite{Ku6}).
The general problem was solved (also in the negative) by Strebel \cite{St2} and the author \cite{Kr6}.
Moreover, due to Strebel, if the initial extremal map $f_0$ has no substantional (essential) boundary points, then, for each $n \ge 4$, equality (38) is attained on $n$-gons only when this map $f_0$ is itself a polygonal map
for some $n$ vertices.

One of the basic tools here is the fundamental Reich-Strebel inequality for polygons, which also will be applied below.

\bigskip\noindent
{\bf 9.2}. Surprisingly, the quadrilaterals and pentagons arising in the cases $n = 4$ and $n = 5$ are intrinsically connected with features of Grunsky norm.
For quadrilaterals ($n = 4$), this was implicitly used by Reich \cite{Re} in a special case.
We shell deal with pentagons; in the case of quadrilaterals the proof is simpler.

Every quasisymmetric homeomorphism of the unit circle $\mathbb S^1$ generates naturally
a quasiconformal automorphism $f^\mu(z)$ conformal on $\D^*$, with dilatation $\|\mu\|_\iy = k_0(h)$
and subject to hydrodynamical normalization $f^\mu(z) = z + b_0 + b_1 z^{-1} + \dots$.

The following theorem connects the extremal pentagonal quasiconformal maps of asymptotically conformal
curves with intrinsic features of Grunsky norm. We fix on $\mathbb S^1$ the points $z_1, \dots, z_5$ regarding the disk $\D$ as polygon $(\D, z_1, \dots, z_n)$ and the function $f^\mu$ as he map of this polygon onto
$f^\mu(\D), f^\mu(z_1, \dots, f^\mu(z_n))$.

\bigskip\noindent
{\bf Theorem 16}. {\it Let a quasicircle $L$ be determined by a quasisymmetric homeomorphism $h$ admitting Teichm\"{u}ller extension to $\D$. Then the equality (38) for $n = 5$ is valid
if and only if its outer conformal map $f: \ \D^* \to D^*$ satisfies
 \be\label{40}
\sup k(f_{[5]}) = \vk(f) = q_L = k_0(h),
\end{equation}
and the extremal extension $f^{\mu_0}$ of $f$ is itself a polygonal map for some $n$ vertices.
}

\bigskip
In particular, all this is valid for asymptotically conformal curves.

\bigskip\noindent
{\bf Proof}. First consider the quasisymmetric automorphisms  $h$ of the real axis $\R$ and fix a pentagon $U(0, 1, a_1, a_2, \iy)$ with the vertices $0, 1, a, b, \iy \in \hR$ such that
$1 < a_1 < a_2 < \iy$.
One can assume that the maps
$f : \ U(0, 1, a, b), \iy)) \to U(h(0), h(1),h(a), h(b), h(\iy))$ fix the points
$0, 1, \iy$.

For any pentagon $U(0, 1, x_1, x_2, \iy)$, the extremal map
$$
f_0 : \ U(0, 1, a, b, \iy)) \to U(0, 1 ,h(a), h(b), \iy)
$$
has the Beltrami coefficient of the form $\mu_0(z) = k |\vp_0(z)|/\vp_0(z)$, where
$$
\vp_0(z) = \fc{c_1}{z(z - 1)(z - a_1)} + \fc{c_2}{z(z - 1)(z - a_2)}
= \frac{(c_1 + c_2) z - (c_1 b + c_2 a)}{z(z - 1)(z - a)(z - b)}
$$
Here the numbers $c_1, c_2$ are real, because the quadratic differential $\vp_0(z) dz^2$ takes
the real values on $\R$, and moreover, $\vp_0$ has only one real zero.

This yields that the extremal quasiconformal maps of any pentagon $D(a_1, a_2, a_3, a_4, a_5)$,
whose domain $D$ is bounded by a $C^{1+\a}$ smooth curve onto other ones, are determind by
holomorphic quadratic differentials
$\psi_0 = (\vp_0 \circ g) (g^\prime)^2$,
where $g$ maps conformally $D(a_1, a_2, a_3, a_4, a_5)$ onto
$U(0, 1, a, b, \infty)$ (or, equvalently, onto the unit disk), and each $\psi_0$ has only one zero; moreover, this zero is located  on the boundary of $D$. Consequently, $\psi_0$ does
not vanish in $D$, and it is possible to define a single-valued branch of $\sqrt{\psi_0}$
in the whole domain $D$.

Now, assuming that equality (38) is valid for a given $h$, one can select a sequence of pentagons
$$
P_m = \D(z_{1,m}, \dots z_{n,m}), \quad m = 1, 2, ... \ ,
$$
and the corresponding extremal polygonal maps
$f_m : P_m \to \D(z_{1,m}, \dots z_{n,m})$
which carry the vertices $z_j$ to $h(z_j), \ j = 1, ... , n$, so that
 \be\label{41}
\lim\limits_{m \to \iy} k(f_m) = k(h) = k(f_0).
\end{equation}
The Beltrami coefficients of these maps are
$\mu_{f_m}(z) = k_m |\vp_m(z)| / \vp_m(z)$ with $k_m = k(f_m)$ and $\vp_m \in A_1(\D), \
\vp_m(z) \ne 0$ in $\D$.
We normalize $\vp_m$ letting $\|\vp_m\|_{A_1(\D)} = 1$.

We now apply to the extremal Beltrami coefficient $\mu_0 = \mu_{f_0}$ of $h$ and to quadratic differentials $\vp_m(z) \ne 0$ in $\D$ the Reich-Strebel inequality for polygons \cite{RS} which yields
$$
\Re \iint_\D
\fc{\mu_0(z) \vp_m(z)}{1 - |\mu_0(z)|^2} \ dx dy \ge \fc{k_m}{1 - k_m}
- \int_\D |\vp_m(z)| \ \fc{|\mu_0(z)|^2}{1 - |\mu_0(z)|^2} \ dx dy.
$$
For the maps $f_0$ with constant $|\mu_0(z)|$ in $U$ (equal to $k(h)$),
the last inequality implies
$$
\fc{1}{1 - k(h)^2} \ \Re \iint_\D \mu_0(z) \vp_m(z) dx dy \ \ge \fc{k_m}{1 - k_m}
- \fc{k(h)^2}{1 - k(h)^2}.
$$
Combining this with (41), one obtains
$$
\liminf\limits_{m \to \iy} \ \Re \iint_\D \mu_0(z) \vp_m(z) dx dy \ge k(h),
$$
which is possible only if
$$
\lim\limits_{m \to \infty} \ \Re \iint_\D \mu_0(z) \vp_m(z) dx dy = k(h).
$$
Since $\vp_m(z) \ne 0$ in $\D$, the last equality implies
 \be\label{42}
\sup\limits_{\vp \in A_1^2(\D) :  \|\vp \|_{A_1(\D)} = 1} \
\Big\vert \iint_\D \mu_0(z) \vp(z) dx dy \Big\vert = k(h),
\end{equation}
and the equalities (40) follow from Theorem 1.

The assertion of theorem that the extremal extension $f^{\mu_0}$ of $f$ must itself be a polygonal map is a consequence of Strebel's result mentioned above.

Conversely, if $\vk(f) = k(f)$ and $f$ has no substantial points on $\mathbb S^1$, this function has
Teichm\"{u}ller extremal extension $f^{\mu_0}$, which obeys (38).

In the case $n = 4$ all terms in (40) are also equal to the quantity (39) defined by conformal moduli of quadrilaterals.

\bigskip\noindent
{\bf 9.3. Strengthened theorem}. The assumption that the extremal extension  of quasisymmetric function  must be or Teichm\"{u}ller type can be omitted, but in the general case the assertion that this
extension is itself an extremal polygonal map is not true.

\bigskip\noindent
{\bf Theorem 17}. {\it The relations (38) and (40)are equivalent for all extremal extensions
of any quasisymmetric function} $h$.

\bigskip\noindent
{\bf Proof}. Let $f^{\mu_0}$ be an extremal extension of a given quasisymmetric homeomorphism
$h$ of $\mathbb S^1$. Since the set of Strebel points is open and dense on the universal
Teichm\"{u}ller space $\T$ (and hence, in the set $\Sigma$), there is a sequence of extremal Teichm\"{u}ller maps $f^{\mu_m}$ satisfying (38) and (40) and such that their Schwarzians
$$
\|S_{f^{\mu_m}} - S_{f^{\mu_0}}\|_\B \to 0 \quad \text{as} \ \ m \to \iy.
$$
The continuity of the Teichm\"{u}ller and Grunsky norms in the distances on $\T$ implies
that such convergence preserves the validity of the indicated relations also for the limit function.
The theorem follows.

\bigskip\noindent
{\bf 9.4. Two consequences}.
First, we mention the following interesting consequence from the proof of Theorem 16.

\bigskip\noindent
{\bf Corollary 2}. {\it The univalent functions $f \in \Sigma$ admitting Teichm\"{u}ller extremal polygonal quasiconformal extensions to $\D$ (for $n = 4, 5$) have equal Teichm\"{u}ller and Grunsky norms, and vise versa.
}

\bigskip
Note that this assertion does not assume that the image $f(\D^*)$ is a standard (unbounded)
polygon; for such polygons the extremal quasiconformal extensions (with the smallest dilatation) can
have the substantional boundary points at vertices and therefore do not be of Teichm\"{u}ller type.

The equality (42) and the known features of functions with equal Grunsky and Teichm\"{u}ller norms
imply

\bigskip\noindent
{\bf Theorem 18}. {\it The strict inequality $k(h) > \sup k(f_{[n]}, n = 4, 5$ is valid on an open
dense subset in the space of quasisymmetric functions.   }

\bigskip\noindent
{\bf 9.5}. In the case $n \ge 6$, the assertions of Theorems 16, 17 are valid only for some
smaller collections extremal Beltrami coefficients (defined by quadratic differentials with zeros of even order). This is illustrated by the following examples.

\bigskip\noindent
{\bf 1}. Pick the polygons  $(\D, e^{i \theta_1}, \dots,  e^{i \theta_n})$ and consider
the rational functions
$$
r_m(z) = p_m(z)/q_m(z),
$$
where $p_m$ and $q_m$ are polynomials of order $m \in \mathbb N$ with zeros placed outside of $D$, and let all zeros of $q_n$ on the unit circle be simple.
Put $\mu_m(z) = |r_m(z)|/r_m(z)$ on $\D$ and extend it by zero to $\D^*$.

This Beltrami coefficient is generated by quasisymmetric homeomorphism $h(\theta)$ of $\mathbb S^1$,
for which the equality (38) is valid.

\bigskip\noindent
{\bf 2}. Consider in $L_1(\D)$ the subspace $\mathcal L_1^0$ spanned by functions
$$
\psi_m(z) = \fc{1}{m} \ e^{- \chi(z)/m}, \quad m = 0, 1, \dots ,
$$
where $\chi$ maps conformally the disk $\D$ onto the right half-strip $\Pi_{+}$ considered in
\textbf{5.3}, and construct on $\mathcal L_1^0$ a linear finctional $l(\psi)$ taking on $\psi_m$
the prescribed values $\eta_m$ with $|\eta_m| < 1, \ m = 0, 1, \dots$, such that $\sum_m |\eta_m| = 1$.
$$
l\Bigl(\sum\limits_0^\iy c_m \psi_m \bigr) = \sum\limits_0^\iy \eta_m \psi_m.
$$
Extending this functional by Hahn-Banach to the entire space $L_1(\D)$, one obtains a bounded
measurable function $\mu$ on $\D$ representing $l$, which satisfies
$$
\|\mu\|_\iy = \sup_m \fc{1}{m} \Big\vert \iint_\D \mu(z) e^{- \chi(z)/m} dxdy \Big\vert.
$$
Now extend $\mu$ by zero to $\D^*$ and take the quasisymmetric homeomorphism o$h(\theta)$ of the unit circle determined by quasiconformal map $f^\mu$. Both functions $h$ and $f^\mu$ satisfy (38).

\bigskip\noindent
{\bf 3}. Let $\mu(z)= r |\psi(z)|/\psi(z)$ on $\D$ with $|r| < 1$ and $\psi \in A_1(\D)$ having
in $\D$ zeros of odd order. Letting
$$
\wt \mu_r(z) = \begin{cases}  t \mu(z),  \ \                          &|z| < 1,   \\
                             \ov{r\mu(1/\ov z)} z^2/\ov z^2 ,     & |z| > 1,
\end{cases}
$$
one obtains $|r|$-quasiconformal automorphism $f^{\wt \mu_r}$ of the sphere $\hC$ preserving
the disks $\D$ and $\D^*$ whose restriction to the unit circle $f^{\wt \mu_r}|\mathbb S^1 = h(t)$ is a quasisymmetric function with
$$
\sup k(h_{[n]}) < k(h) = k(f^{\wt \mu_r}) = |r|.
$$

The special case $\psi(z) = z$ was considered in \cite{Re}.

\bigskip\bigskip
\centerline{\bf 10. OPEN PROBLEMS}

\bigskip\noindent
{\bf 1}. Extend the Pommerenke-Zhuravlev theorem and Theorem 9 to univalent functions in arbitrary quasidisks.

\bigskip\noindent
{\bf 2}. Does the equality (38) be valid (for any $n \ge 4)$ for the boundary values of all univalent functions $f(z)$ with $\vk(f) = k(f)$ ?

\bigskip\noindent
{\bf 3}. The above theory is intrinsically connected with the standard classes of univalent functions
having quasiconformal extensions. How about its connection with other and more general classes of qusiconformal maps, for example, with extremality of solutions of quasilinear Beltrami equations considered by V. Gutlyanskii and V. Ryazanov in \cite{GR}?

\bigskip
{\small\em{ \leftline{Department of Mathematics, Bar-Ilan
University, 5290002 Ramat-Gan, Israel} \leftline{and
Department of Mathematics, University of Virginia,  Charlottesville, VA 22904-4137, USA}}

\end{document}